\documentclass[11pt]{article}

\usepackage{bm}
\usepackage{epsfig}
\usepackage{caption}
\usepackage{amsmath,amsfonts} 
\usepackage{graphics}                 
\usepackage{graphicx}
\usepackage{float}
\usepackage{array}
\usepackage{epsf}
\usepackage{color}                
\usepackage{accents}

\newtheorem{theorem}{Theorem}[section]
\newtheorem{lemma}[theorem]{Lemma}
\newtheorem{remark}[theorem]{Remark}

\newtheorem{example}{Example}[section]

\newenvironment{proof}[1][Proof]{\begin{trivlist}
\item[\hskip \labelsep {\bfseries #1}]}{\end{trivlist}}

\newcommand{\nn}{\nonumber}

\def\refe#1{(\ref{#1})}
\def \bb#1{\setbox0=\hbox{$#1$}
\kern-.025em\copy0\kern-\wd0 \kern.05em\copy0\kern-\wd0
\kern-.025em\raise.0433em\box0}
\def \endproof{\vrule height8pt width 5pt depth 0pt}

\newcommand{\qed}{\nobreak \ifvmode \relax \else
      \ifdim\lastskip<1.5em \hskip-\lastskip
      \hskip1.5em plus0em minus0.5em \fi \nobreak
      \vrule height0.75em width0.5em depth0.25em\fi}

\newcommand{\bA}{\mathbf{A}}
\newcommand{\bv}{\mathbf{v}}
\newcommand{\bH}{\mathbf{H}}
\newcommand{\D}{\mathrm{div}}
\newcommand{\C}{\mathbf{curl}}
\newcommand{\curl}{\mathrm{curl}}

\newcommand{\n}{\mathbf{n}}

\newcommand{\ee}{\epsilon}

\begin{document}

\title{Analysis of linearized Galerkin-mixed FEMs
for the time-dependent Ginzburg--Landau equations 
of superconductivity}

\author{
\setcounter{footnote}{0}
Huadong~Gao
\footnote{
School of Mathematics and Statistics, 
Huazhong University of Science and Technology, 
Wuhan 430074, P.R. China.
{\tt huadong@hust.edu.cn}.
The work of the author was supported in part by the 
National Science Foundation of China No. 11501227.}
\quad and \quad Weiwei~Sun 
\footnote{ 
Department of Mathematics,
City University of Hong Kong, Kowloon, Hong Kong.
{\tt maweiw@math.cityu.edu.hk}.
The work of the author was supported in part by a grant
from the Research Grants Council of the Hong Kong Special
Administrative Region, China. (Project No. CityU 11302514).}
}
 
\date{}

\maketitle

\begin{abstract}
A linearized backward Euler Galerkin-mixed finite element method
is investigated for the time-dependent Ginzburg--Landau (TDGL) 
equations under the Lorentz gauge.
By introducing the induced magnetic field $\bm{\sigma} = \C \, \bA$
as a new variable, the Galerkin-mixed FE scheme 
offers many advantages over conventional Lagrange
type Galerkin FEMs. An optimal error estimate for
the linearized Galerkin-mixed FE scheme is established unconditionally.
Analysis is given under more general assumptions for the regularity of
the solution of the TDGL equations, which includes the problem in
two-dimensional noncovex polygons and certain three dimensional polyhedrons, 
while the conventional Galerkin FEMs may not converge to a true solution in these cases.
Numerical examples in both two and three dimensional spaces are
presented to confirm our theoretical analysis.
Numerical results show clearly the efficiency of the mixed method,
particularly for problems on nonconvex domains.

\vskip 0.2in
\noindent{\bf Keywords:} 
Ginzburg--Landau equation, linearized scheme, mixed finite element method,
unconditional convergence, optimal error estimate, superconductivity.

\end{abstract} 

\section{Introduction}
\setcounter{equation}{0}
In this paper, 
we consider the time-dependent Ginzburg--Landau (TDGL) 
equations under the Lorentz gauge
\begin{eqnarray}
&&\eta \frac{\partial \psi}{\partial t} 
- i \eta \kappa ({\D \, \bA}) \psi 
+ (\frac{i}{\kappa} \nabla + \bA)^{2} \psi 
+ (|\psi|^{2}-1) \psi = 0 \, ,
\label{pde1}\\
&&\frac{\partial \bA}{\partial t} - \nabla {\D \, \bA} 
+ {\C \, \C \, \bA} 
+ \frac{i}{2 \kappa}(\psi^{*} \nabla \psi - \psi \nabla \psi^{*}) 
+ |\psi|^{2} \bA = \C \, \mathbf{H}_e \, , 
\label{pde2}
\end{eqnarray}
for $x \in \Omega$ and $t \in (0,T]$, where 
the complex scalar function $\psi$ is the order parameter
and the real vector-valued function $\bA$ is the magnetic potential.
In \refe{pde1}-\refe{pde2}, $|\psi|^2$ denotes the density of the
superconducting electron pairs. $|\psi|^2 = 1$ and $|\psi|^2 = 0$
represent the perfectly superconducting state and
the normal state, respectively, while $ 0 < |\psi|^2 < 1$ represents a mixed state.
The real vector-valued function $\bH_e$ is the external applied magnetic field,
$\kappa$ is the Ginzburg--Landau parameter and 
$\eta$ is a dimensionless constant. In the following, we set $\eta = 1$ for
the sake of simplicity.
We assume that $\Omega$ is a simply-connected bounded Lipschitz domain 
in $\mathbb{R}^{3}$.
The following boundary and initial conditions are 
supplemented to \refe{pde1}-\refe{pde2}
\begin{align}
&\frac{\partial \psi}{\partial \n} = 0, \quad 
\C \, \bA \times \mathbf{n} = \mathbf{H}_e \times \mathbf{n} , \quad
\bA \cdot \mathbf{n} = 0, 
&&\mathrm{on}\,  \partial \Omega \times [0,T],
\label{original-bc}\\
&\psi(x,0) = \psi_{0}(x), \quad \bA(x,0) = \bA_{0}(x), \quad 
&& \mathrm{in} \, \Omega,
\label{original-init}
\end{align}
where $\mathbf{n}$ is the outward unit normal vector.
We note that in \refe{pde1}-\refe{original-init}, 
$\nabla$, $\D$ and $\C$ are the standard gradient, divergence
and rotation operators in three dimensional space.

The TDGL equations were first deduced by Gor'kov and Eliashberg 
in {\cite{GE}} from the microscopic Bardeen--Cooper--Schrieffer theory
of superconductivity. 
For detailed physical description and mathematical modeling
of the superconductivity phenomena, we refer to the review articles 
{\cite{CHO,du}}. 
Theoretical analysis for the TDGL equations 
can be found in literature \cite{CHL,Feireisl-Takac,LO,RWW,TW}.
The global existence and uniqueness of 
the strong solution were established in  \cite{CHL} 
for the TDGL equations with the Lorentz gauge. 
Numerical methods for solving the TDGL equations  
have also been studied extensively, e.g., see  
\cite{AlstromSPM,CH,Du_first,Dureview,GLS,gropp,GunKL,Mu,MH,RPCA,win}. 
A semi-implicit Euler scheme with a finite element 
approximation was first proposed by Chen and Hoffmann \cite{CH} 
for the TDGL equations with the Lorentz gauge. 
A suboptimal $L^2$ error estimate  
was obtained for the equations in two dimensional space. 
Later, a decoupled alternating Crank--Nicolson Galerkin method was 
proposed by Mu and Huang {\cite{MH}}. An optimal error estimate was 
presented under the time step restrictive conditions 
$\tau = O(h^{\frac{11}{12}})$ for the two-dimensional model and 
$\tau = O(h^{2})$ for the three-dimensional model, 
where $h$ and $\tau$ are the mesh size 
in the spatial direction and the time direction, respectively. 
All these schemes in {\cite{CH,Du_first,MH}}
are nonlinear.  At each time step, one has to solve a nonlinear system.  
Clearly, more commonly-used discretizations
for nonlinear parabolic equations are linearized schemes, 
which only require the solution of a linear system at each time step.  
A linearized Crank--Nicolson type scheme was proposed in \cite{Mu} 
for slightly different TDGL equations 
and a systematic numerical investigation was made there. 
An optimal error estimate of a linearized Crank--Nicolson Galerkin FE scheme
for the TDGL equations was provided unconditionally in a recent work 
\cite{GLS}.

All these methods mentioned above were introduced to solve the TDGL
equations \refe{pde1}-\refe{original-init}  for the order parameter $\psi$ and
the magnetic potential $\bA$ by Galerkin FEMs (or finite
difference methods) and then, to calculate
the magnetic field $\C \, \bA$ by certain numerical differentiation.
There are several drawbacks for such approaches. 
One of important issues in the study of the vortex motion 
in superconductors is the influence of  geometric defects, 
which is of high interest in physics and can be considered as a problem in
polygons. 
Previous numerical results \cite{GS, GLS, Mu} by conventional
Lagrange FEMs showed the singularity and inaccuracy of
the numerical solution around corners of polygons.
Moreover, for the problem in a nonconvex polygon, 
conventional Lagrange FEMs for the TDGL equations may
converge to a spurious solution, 
see the numerical experiments reported in \cite{GS,Li_Zhang2} for
the problem on an $L$-shape domain.
Another issue in the superconductivity model is
its coupled boundary conditions in \refe{original-bc},
which may bring some extra difficulties in implementation to conventional
Lagrange FEMs on a general domain,
see section 5 of \cite{chen3} for more detailed description.
It is natural that a mixed finite element method with
the physical quantity $\bm{\sigma} = \C \, \bA$ being the new variable may
offer many advantages. 

It has been noted that mixed methods for second order elliptic problems, 
such as the scalar Poisson equation $-\Delta u  =f$ with 
simple boundary conditions, Dirichlet or Neumann, 
have been well studied with the extra variable 
$\bm{v}=\nabla u$ in the last several decades, 
e.g., see \cite{BS,BF,Gat,GR,thomee} and references therein. 
However, mixed finite element methods for the vector 
Poisson equation $-\Delta {\bf u} = {\bf f}$ 
with the coupled boundary conditions
${\bf curl }\, {\bf u} \times {\bf n} = {\bf g}$ 
and ${\bf u} \cdot \mathbf{n} = 0$
seems much more complicated than for a scalar equation. 
Theoretical analysis was established quite recently \cite{AFW,AFW1}
in terms of finite element exterior calculus. 
For the TDGL equations \refe{pde1-d2}-\refe{init-d2} in two dimensional space, 
Chen \cite{chen3} proposed a semi-implicit and weakly nonlinear 
scheme with a mixed finite element method, in which 
both $\curl \, \bA$ and $\D \, \bA$ were introduced to be
extra unknowns, and the corresponding finite element spaces for 
$\curl \, \bA$ and $\D \, \bA$ were constructed with certain bubble functions
and some special treatments in elements adjacent to the boundary of the domain 
were also needed. 
A suboptimal $L^2$ error estimate was provided 
and no numerical example was given in \cite{chen3}. 
Recently, a linearized backward Euler Galerkin-mixed
finite element method was proposed in \cite{GS} 
for the TDGL equations in both two
and three dimensional spaces with only one extra physical quantity
$\bm{\sigma} = \C \, \bA$ (${\sigma} = \curl \, \bA$ in two dimensional space). 
No analysis was provided in \cite{GS}. 
The scheme is decoupled and at each time step,
one only needs to solve two linear systems for $\psi$ and
$(\bA, \bm{\sigma})$, which can be done simultaneously.

This paper focuses on theoretical analysis on
a linearized Galerkin-mixed FEM for the TDGL equations to
establish optimal error estimates.
More important is that our analysis is given
without any time-step restrictions and 
under more general assumptions for the regularity of
the solution of the TDGL equations, which includes the problem in
nonconvex polygons and certain convex polyhedrons. 
Usually, conventional FEMs for a scalar parabolic equation
require the regularity of the solution in $H^{1+s}$ with 
$s>0$ \cite{CLTW,CHou}, 
while for the TDGL equations in a nonconvex polygon, 
$\bA \in \bH^s$ and $\curl \, \bA$, $\D \, \bA \in H^{1+s}$ with $s<1$ 
in general \cite{Li_Zhang}. 
Our numerical results show clearly that the mixed method converges to
a true solution for problems in a nonconvex polygon and 
conventional Lagrange type FEMs do not in this case. 
In addition, in the mixed FEM, the coupled boundary condition
$\C \, \bA \times \mathbf{n} = \mathbf{H}_e \times \mathbf{n}$
reduces to a Dirichlet type boundary condition
$\bm{\sigma} \times \n = \mathbf{H}_e \times \n$.
Implementation of Dirichlet boundary conditions for vector-valued elements
(such as {N\'ed\'elec} and Raviart--Thomas elements) becomes much simpler.
In a recent work by Li and Zhang, an approach based on Hodge decomposition 
was proposed and analyzed with optimal error estimates. 
However, the approach is applicable only for 
the problem in two-dimensional space due to the restriction of the Hodge 
decomposition. 

The rest of this paper is organized as follows. 
In section \ref{femmethod}, we introduce 
the linearized backward Euler scheme 
with Galerkin-mixed finite element approximations for the TDGL equations 
and we present our main results on unconditionally optimal error estimates. 
In section \ref{preparation}, we recall some results for 
an auxiliary elliptic problem of the vector Poisson equation 
with the coupled boundary conditions.
In section \ref{proofmainresult}, 
we prove that an optimal $L^2$ error estimate
holds almost unconditionally ($i.e.$, without any mesh ratio restriction). 
In section \ref{numericalresults}, 
we provide several numerical examples to confirm our theoretical analysis 
and show the efficiency of the proposed methods. 
Some concluding remarks are given in section \ref{sec-conclusion}.

\section{ A linearized backward Euler Galerkin-mixed FEM}
\label{femmethod}
\setcounter{equation}{0}

In this section, we present a linearized backward Euler
Galerkin-mixed finite element method for the TDGL equations
and our main results. 
For simplicity, we introduce some standard notations and operators below. 
For any two complex functions $u$, $v \in \mathcal{L}^{2}(\Omega)$, 
we denote the $\mathcal{L}^{2}(\Omega)$ inner product and norm by
\begin{equation}
(u,v)  = \int_{\Omega} u(x) \, (v(x))^{*} \, {\mathrm{d}} x, 
\qquad {\left \| u \right \|_{L^2}} = (u,u)^{\frac{1}{2}} \, ,
\nn
\end{equation}
where $v^*$ denotes the conjugate of the complex function $v$. 
Let $W^{k,p}(\Omega)$ be the Sobolev space defined on $\Omega$, 
and by conventional notations, $H^{k}(\Omega) := W^{k,2}(\Omega)$,
$\accentset{\circ}{H}^{k}(\Omega) := \accentset{\circ}{W}^{k,2}(\Omega)$.
Let $\mathcal{H}^{k}(\Omega)=\{u+iv| \, u,v \in H^{k}(\Omega)\}$ be
a complex-valued Sobolev space
and $\bH^{k}(\Omega) = [H^{k}(\Omega)]^{d}$ 
be a vector-valued Sobolev space,
where $d$ is the dimension of $\Omega$.
For a positive real number $s = k + \theta$ with $0 < \theta < 1$,
we define $H^s(\Omega)= (H^k, H^{k+1})_{[\theta]}$ by the complex interpolation,
see \cite{Bergh_Lofstrom}.
To introduce the mixed variational formulation,
we denote 
\begin{align}
{\bH}(\D) = \left\{
\bA \, \big| \, \bA \in \mathbf{L}^{2}(\Omega),
\D \, \bA  \in L^2(\Omega) \right\} \,
\textrm{with $\|\bA\|_{{\bH}(\D)} = \left(\|\bA\|_{{L}^2}^2 
+ \|\D \, \bA\|_{{L}^2}^2  \right)^{\frac{1}{2}} $}
\nn
\end{align}
and 
\begin{align}
{\bH}(\C) = \left\{
\bA \, \big| \, \bA \in \mathbf{L}^{2}(\Omega),
\C \, \bA  \in \mathbf{L}^2(\Omega) \right\} \,
\textrm{with $\|\bA\|_{{\bH}(\C)} = \left(\|\bA\|_{{L}^2}^2 
+ \|\C \, \bA\|_{{L}^2}^2  \right)^{\frac{1}{2}} $} .
\nn
\end{align}
Also we define 
\begin{align}
\accentset{\circ}{\bH}(\D) = \left\{
\bA \, \big| \, \bA \in {\bH}(\D), \quad
\bA \cdot \mathbf{n} \big|_{\partial \Omega} = 0
\right\} 
\nn
\end{align}
and its dual space $\accentset{\circ}{\bH}(\D)'$
with norm 
\begin{align}
\left \| \mathbf{v} \right \|_{\accentset{\circ}{\bH}(\D)'}
:= \sup_{\mathbf{w} \in \accentset{\circ}{\bH}(\D)}
\frac{(\mathbf{v} \, , \, \mathbf{w})}
{\left \| \mathbf{w} \right \|_{{\bH}(\D)}}\,.
\nn
\end{align}
Moreover, we denote 
\begin{align}
\accentset{\circ}{\bH}(\C) = \left\{
\bA \, \big| \, \bA \in {\bH}(\C), \quad
\bA \, \times \, \mathbf{n} \big|_{\partial \Omega} = \mathbf{0} 
\right\}\,.
\nn
\end{align}

By introducing $\bm{\sigma} = \C \,\bA$, the mixed form
of the TDGL equations \refe{pde1}-\refe{original-init}
can be written by
\begin{eqnarray}
&& \frac{\partial \psi}{\partial t} 
- i \kappa ({\D \, \bA}) \psi 
+ (\frac{i}{\kappa} \nabla + \bA)^{2} \psi 
+ (|\psi|^{2}-1) \psi = 0 \, ,
\label{mixed-pde1}\\[6pt]
&& \bm{\sigma} = \C \bA \, ,
\label{mixed-pde1.5} \\
&&\frac{\partial \bA}{\partial t} - \nabla {\D \, \bA} 
+ {\C \, \bm{\sigma}} 
+ \frac{i}{2 \kappa}(\psi^{*} \nabla \psi - \psi \nabla \psi^{*}) 
+ |\psi|^{2} \bA = \C \, \mathbf{H}_e \, ,
\label{mixed-pde2}
\end{eqnarray}
with boundary and initial conditions
\begin{align}
& \frac{\partial \psi}{\partial \n} = 0, \quad 
\bm{\sigma} \times \n = \mathbf{H}_e \times \n, \quad
\bA \cdot \mathbf{n} = 0, \quad 
&& \textrm{on $\partial \Omega \times [0,T]$},
\label{bc}\\
& \psi(x,0) = \psi_{0}(x), \quad
\bm{\sigma} (x,0) = \C \, \bA_{0}(x), \quad
\bA(x,0) = \bA_{0}(x), \quad
&& \mathrm{in}\ \Omega.
\label{init}
\end{align}
The mixed variational formulation of 
the TDGL equations (\ref{mixed-pde1})-(\ref{mixed-pde2}) 
with boundary and initial conditions 
\refe{bc}-\refe{init}
is to find $\psi \in L^2(0,T;\mathcal{H}^{1}(\Omega))$
with $\frac{\partial \psi}{\partial t} \in L^2(0,T;\mathcal{H}^{-1}(\Omega))$,
and $(\bm{\sigma}, \bA) \in L^2(0,T;\bH(\C)) \times 
L^2(0,T;\accentset{\circ}{\bH}(\D)) $ with
$\frac{\partial \bA}{\partial t} \in L^2(0,T;\accentset{\circ}{\bH}(\D)')$,
where $\bm{\sigma} \times \n = \mathbf{H}_e \times \n $ on $\partial \Omega$,
such that 
\begin{align} 
(\frac{\partial \psi}{\partial t}, \omega)
& - i\kappa (({\D \, \bA}) \psi,\omega) 
+ ((\frac{i}{\kappa} \nabla + \bA) \psi,(\frac{i}{\kappa} \nabla 
+ \bA) \omega) 
\nn \\
& + ((|\psi|^{2}-1) \psi, \omega) = 0 \, , 
\quad \forall  \, \omega \in \mathcal{H}^{1}(\Omega)
\label{variation1}
\end{align} 
and
\begin{align} 
& (\bm{\sigma}, \bm{\chi}) - (\C \, \bm{\chi}, \bA) = 0, \quad
  \forall \, \bm{\chi} \in \accentset{\circ}{\bH}(\C)\, ,
\label{variation1.5}
\\[5pt]
&(\frac{\partial \bA}{\partial t},\bv) 
+ ( \C \, \bm{\sigma}, \bv) 
+ ({\D \, \bA},{\D \, \bv}) 
+ \frac{i}{2 \kappa}((\psi^{*} \nabla \psi - \psi \nabla \psi^{*}), \bv)
\nn \\
& \qquad \qquad + (|\psi|^{2} \bA,\bv) 
= (\C \, \mathbf{H}_e, \bv) \, ,  \quad 
\forall \, \bv \in \accentset{\circ}{\bH}(\D) \, ,
\label{variation2}
\end{align} 
for a.e. $t \in (0,T]$
with $\psi(x,0) = \psi_{0}(x)$, $\bm{\sigma} (x,0) = \C \, \bA_{0}(x)$ 
and $\bA(x,0) = \bA_{0}(x)$.

For simplicity, we assume that $\Omega$ is a polyhedron in three 
dimensional space.
Let $\mathcal{T}_{h}$ be a quasi-uniform tetrahedral partition of $\Omega$ 
with $\Omega = \cup_{K} \Omega_{K}$ and denote by  
$h=\max_{\Omega_{K} \in \mathcal{T}_{h}} 
\{ \mathrm{diam} \, \Omega_{K}\}$ the mesh size. 
For a given partition $\mathcal{T}_{h}$, 
we denote by $\mathcal{V}_{h}^{r}$ the $r$-th order Lagrange finite element 
subspace of $\mathcal{H}^{1}(\Omega)$.
we denote by ${\mathbf{Q}}_{h}^{r}$ the $r$-th 
order first type {N\'ed\'elec} finite element subspace of 
$\bH(\C)$, where the case $r=1$ corresponds to 
the lowest order {N\'ed\'elec} edge element ($6$ dofs). 
We denote by ${\mathbf{U}}_{h}^{r}$ the $r$-th 
order Raviart--Thomas finite element subspace of 
${\bH}(\D)$, where the case $r=0$ corresponds to 
the lowest order Raviart--Thomas face element ($4$ dofs).
We also define $\accentset{\circ}{\mathbf{Q}}_{h}^{r} = 
\mathbf{Q}_{h}^{r}  \cap \accentset{\circ}{\bH}(\C)$
and $\accentset{\circ}{\mathbf{U}}_{h}^{r} = 
\mathbf{U}_{h}^{r}  \cap \accentset{\circ}{\bH}(\D)$.
It should be remarked that Brezzi--Douglas--Marini element can also
be used for approximation of ${\bH}(\D)$. 
Here we only confine our attention to 
the Raviart--Thomas element approximation.
By noting the approximation properties of the finite element spaces
$\mathcal{V}_{h}^{r}$, $\mathbf{Q}_{h}^{r}$ and 
$\mathbf{U}_{h}^{r}$ \cite{Alonso_Valli,BS,Gat},
we denote by $\pi_h$ a general projection operator on  $\mathcal{V}_{h}^{r}$, 
$\mathbf{Q}_{h}^{r}$ and $\mathbf{U}_{h}^{r}$, satisfying 
\begin{equation}
\left \{
\begin{array}{ll}
{\left\|\omega  - \pi_h \omega \right\|}_{L^{2}}
\le C h^s {\left\| \omega \right\|}_{H^{s}} \, , 
& 0 < s \le r+1 \, ,
\\[4pt]
{\left\|\bm{\chi}  - \pi_h \bm{\chi}\right\|}_{L^{2}} +
{\left\|\C \, (\bm{\chi} - \pi_h \bm{\chi})\right\|}_{L^{2}} 
\le C h^s ({\left\|\bm{\chi}\right\|}_{H^{s}} +
{\left\|\C \, \bm{\chi}\right\|}_{H^{s}}), 
& \frac{1}{2} < s \le r \, ,
 \\[4pt]
{\left\|\bm{v}  - \pi_h \bm{v}\right\|}_{L^{2}} +
{\left\|\D \, (\bm{v} - \pi_h \bm{v})\right\|}_{L^{2}} 
\le C h^s ({\left\|\bm{v}\right\|}_{H^{s}} +
{\left\|\D \, \bm{v}\right\|}_{H^{s}}), 
& 0 < s \le r+1 \, . 
\end{array}
\right.
\label{interpolate}
\end{equation}

Let ${\left \{ t_{n} \right \}}_{n=0}^{N}$ be a uniform partition
in the time direction with the step size $\tau = \frac{T}{N}$,
and let $u^n = u(\cdot, n \tau)$. 
For a sequence of functions $\{U^{n}\}_{n=0}^{N}$ defined on $\Omega$, 
we denote 
\begin{eqnarray}
{D_{\tau}} U^{n} = \frac{U^{n}-U^{n-1}}{\tau}, \quad
\textrm{for $n=1$, $2$, $\ldots$, $N$}.
\nn
\end{eqnarray}

With the above notations, 
the linearized backward Euler Galerkin-mixed FEM
for the mixed form TDGL equations \refe{mixed-pde1}-\refe{init}
is to find $\psi_{h}^{n} \in \mathcal{V}_{h}^{\widehat{r}}$
and $(\bm{\sigma}_{h}^{n} \, , \bA_{h}^{n}) 
\in \mathbf{Q}_h^{r+1} \times \accentset{\circ}{\mathbf{U}}_{h}^{r}$,
with $\bm{\sigma}_{h}^{n} \times \n = \pi_h\mathbf{H}_e^n \times \n$ on $\partial \Omega$,
such that for $n = 1,2, \ldots, N$, 
\begin{eqnarray}
&&(D_{\tau}\psi_{h}^{n}, \omega_{h}) 
- i\kappa ((\D \,\bA_{h}^{n-1}) \psi_{h}^{n},\omega_{h}) 
+ ( (\frac{i}{\kappa}\nabla+\bA_h^{n-1}) {\psi}_{h}^{n} \,, 
(\frac{i}{\kappa}\nabla+\bA_h^{n-1}) \omega_{h})
\nn \\
&&
\qquad\qquad +((|\psi_{h}^{n-1}|^{2}-1)\psi_{h}^{n}, \omega_{h}) =0 \, ,
\quad \forall \, \omega_{h} \in \mathcal{V}_{h}^{\widehat{r}} \, , 
\label{fem1} 
\end{eqnarray} 
and
\begin{eqnarray} 
&& (\bm{\sigma}_h^n, \bm{\chi}_h) - (\C \, \bm{\chi}_h, \bA_h^n) = 0, 
\quad \forall  \, 
\bm{\chi}_h \in \accentset{\circ}{\mathbf{Q}}^{r+1}_{h}\, , 
\label{fem1.5}
\\[8pt]
&&(D_{\tau}\bA_{h}^{n},\bv_{h}) 
+ ({\D \, {\bA}_{h}^{n}} \, ,{\D \, \bv_{h}})
+ ( \C \, \bm{\sigma}_h^n \, , \bv_{h})
+ (|{\psi}_{h}^{n-1}|^{2} {\bA}_{h}^{n},\bv_{h})
\nn \\
&&  = (\C \, \mathbf{H}_e^{n} \, ,\bv_{h}) 
- \frac{i}{2 \kappa} \left( ({\psi}_{h}^{n-1})^{*} 
\nabla {\psi}_{h}^{n-1} 
- {\psi}_{h}^{n-1} \nabla ({\psi}_{h}^{n-1})^{*}, \bv_{h} \right) \, ,
\quad \forall \, \bv_{h} \in \accentset{\circ}{\mathbf{U}}_{h}^{r} \, ,
\label{fem2}
\end{eqnarray}
where $r \ge 0$ and $\widehat{r} = \max \{ 1,r \}$.
$\psi_{h}^{0} = \pi_h\psi_{0}$ and $\bA_{h}^{0} = \pi_h \bA_0$
are used at the initial time step. 

The linearized backward Euler Galerkin-mixed FEM scheme (\ref{fem1})-(\ref{fem2}) 
is uncoupled and $\psi_h^n$ and $(\bm{\sigma}_h^n \, , \bA_h^n)$
can be solved simultaneously. Moreover, for each FEM equation, one only needs to solve 
a linear system at each time step.

\begin{remark}
{\bf A linearized Galerkin-mixed FEM in two dimensional space}:
If the superconductor is a long cylinder in the $z$-direction
with a finite cross section and the external applied 
field $\mathbf{H}_e = H_e [0,0,1]^T$ (i.e., 
$\mathbf{H}_e$ is parallel to the $z$-axis), the original three dimensional
equations \refe{pde1}-\refe{pde2} can be reduced to 
a two dimensional equation \cite{GunKL,RPCA}
\begin{eqnarray}
&& \frac{\partial \psi}{\partial t} 
- i \kappa ({\D \, \bA}) \psi 
+ (\frac{i}{\kappa} \nabla + \bA)^{2} \psi 
+ (|\psi|^{2}-1) \psi = 0 \, ,
\label{pde1-d2}\\
&&\frac{\partial \bA}{\partial t} - \nabla {\D \, \bA} 
+ {\C \, \curl \, \bA} 
+ \frac{i}{2 \kappa}(\psi^{*} \nabla \psi - \psi \nabla \psi^{*}) 
+ |\psi|^{2} \bA = \C \, H_e \, , 
\label{pde2-d2}
\end{eqnarray}
with boundary and initial conditions 
\begin{eqnarray}
&&\frac{\partial \psi}{\partial \n} = 0, \quad 
\curl \, \bA = H_e, \quad \bA \cdot \mathbf{n} = 0, \quad 
\mathrm{on}\  \partial \Omega \times (0,T],
\label{bc-d2}\\
&&\psi(x,0) = \psi_{0}(x), \quad \bA(x,0) = \bA_{0}(x), \quad 
\mathrm{in}\ \Omega \, ,
\label{init-d2}
\end{eqnarray}
where $\psi$ and $\bA = [A_1,A_2]^T$ are
scalar-valued complex function and vector-valued real function, respectively.
The operators $\D$, $\nabla$, $\curl$ and $\C$ in \refe{pde1-d2}-\refe{init-d2}
are defined by
\begin{align}
&\D \, \bA = \frac{\partial A_1}{\partial x}
+ \frac{\partial A_2}{\partial y} , \, 
\nabla \psi  = 
\left[\frac{\partial \psi}{\partial x}\, ,
\frac{\partial \psi}{\partial y} \right]^{T} , \,
 \curl \, \bA = \frac{\partial A_2}{\partial x}
- \frac{\partial A_1}{\partial y} \, , 
\C \, \psi = \left[\frac{\partial \psi}{\partial y}\, ,
-\frac{\partial \psi}{\partial x} \right]^{T}  .
\nn
\end{align}
To introduce the FEM scheme for \refe{pde1-d2}-\refe{init-d2},
for a quasi-uniform triangular mesh $\mathcal{T}_{h}$ with
mesh size $h=\max_{\Omega_{K} \in \mathcal{T}_{h}} \{ \mathrm{diam} \, \Omega_{K}\}$,
we denote by $\mathcal{V}_{h}^{r}$ and ${\mathbf{Q}}_{h}^{r}$
the $r$-th order Lagrange finite element 
subspaces of $\mathcal{H}^{1}(\Omega)$ and ${H}^{1}(\Omega)$, respectively.
We also define $\accentset{\circ}{\mathbf{Q}}_h^{r} = 
{\mathbf{Q}}_{h}^{r} \cap \accentset{\circ}{H}^{1}(\Omega)$.
Let $\accentset{\circ}{\mathbf{U}}_{h}^{r}$ be the 
$r$-order Raviart--Thomas finite element subspaces of $\accentset{\circ}{\bH}(\D)$,
where the case $r=0$ corresponds to 
the lowest order Raviart--Thomas element ($3$ dofs).
Then, by introducing  ${\sigma} = \curl \,\bA$,
the linearized backward Euler Galerkin-mixed FEM scheme
is to find $\psi_{h}^{n} \in \mathcal{V}_{h}^{\widehat{r}}$
and $({\sigma}_{h}^{n} \, , \bA_{h}^{n}) 
\in \mathbf{Q}_h^{r+1} \times \accentset{\circ}{\mathbf{U}}_{h}^{r}$,
where ${\sigma}_{h}^{n} = \pi_h {H_e^n}$ on $\partial \Omega$,
such that for $n = 1,2, \ldots, N$
\begin{align} 
(D_{\tau}\psi_{h}^{n}, \omega_{h}) 
& - i\kappa ((\D \,\bA_{h}^{n-1}) \psi_{h}^{n},\omega_{h}) 
+ ( (\frac{i}{\kappa}\nabla+\bA_h^{n-1}) {\psi}_{h}^{n} \,, 
(\frac{i}{\kappa}\nabla+\bA_h^{n-1}) \omega_{h})
\nn \\
&
+((|\psi_{h}^{n-1}|^{2}-1)\psi_{h}^{n}, \omega_{h}) =0 \, ,
\quad \forall \, \omega_{h} \in \mathcal{V}_{h}^{\widehat{r}} 
\label{fem1-d2} 
\end{align} 
and
\begin{eqnarray} 
&& ({\sigma}_h^n, {\chi}_h) - (\C \, {\chi}_h, \bA_h^n) = 0, 
\quad \forall  \, 
{\chi}_h \in \accentset{\circ}{\mathbf{Q}}^{r+1}_{h}\, , 
\label{fem1.5-d2}
\\[8pt]
&&(D_{\tau}\bA_{h}^{n},\bv_{h}) 
+ ({\D \, {\bA}_{h}^{n}} \, ,{\D \, \bv_{h}})
+ ( \C \, \sigma_h^n \, , \bv_{h})
+ (|{\psi}_{h}^{n-1}|^{2} {\bA}_{h}^{n},\bv_{h})
\nn \\
&&  = (\C \, {H}_e^{n} \, ,\bv_{h}) 
- \frac{i}{2 \kappa} \left( ({\psi}_{h}^{n-1})^{*} 
\nabla {\psi}_{h}^{n-1} 
- {\psi}_{h}^{n-1} \nabla ({\psi}_{h}^{n-1})^{*}, \bv_{h} \right) \, ,
\quad \forall \, \bv_{h} \in \accentset{\circ}{\mathbf{U}}_{h}^{r} \, , 
\label{fem2-d2}
\end{eqnarray}
where $\psi_{h}^{0} = \pi_h\psi_{0}$ and $\bA_{h}^{0} = \pi_h \bA_0$
are used at the initial time step. 
\end{remark}

Here we focus our attention on analysis of the mixed scheme and 
present our main results on optimal error estimates in the following theorem. 
The proof for the problem in three dimensional space 
will be given in sections 3 and 4. 
The proof for the two-dimensional model can be obtained analogously and 
therefore, omitted here. 
Numerical simulations on a slightly different scheme were given in \cite{GS}. 
Further comparison with conventional Galerkin FEMs will be presented in 
the section \ref{numericalresults}.
We assume that the initial-boundary
value problem (\ref{pde1})-(\ref{original-init}) 
has a unique solution satisfying the regularity
\begin{eqnarray}
&& \psi\in {L^{\infty}(0,T;\mathcal{H}^{{l+1}})} \, ,
\psi_{t} \in {L^{\infty}(0,T;\mathcal{H}^{{l+1}})} \, , 
\psi_{tt} \in {L^{\infty}(0,T;\mathcal{L}^2)}
\label{regularity-1}
\end{eqnarray}
and
\begin{equation} 
\left \{
\begin{array}{l}
\bA \in {L^{\infty}(0,T;\mathbf{H}^{l})} \, , 
\bA_{t} \in {L^{\infty}(0,T;\mathbf{H}^{l})} \, , 
\bA_{tt} \in {L^{\infty}(0,T;\mathbf{L}^2)} \, ,
\\
 \D \, \bA \in {L^{\infty}(0,T;{H}^{l+1})} \, , 
\left(\D \,\bA \right)_{t} \in {L^{\infty}(0,T;{H}^{l+1})} \, ,
\\
\bm{\sigma} \in {L^{\infty}(0,T;\mathbf{H}^{l+1})} \, , 
\bm{\sigma}_{t} \in {L^{\infty}(0,T;\mathbf{H}^{l+1})} \, , 
\label{regularity-2}
\end{array}
\right.
\end{equation} 
where $l > \frac{1}{2}$ depends on the regularity of the domain $\Omega$.

\begin{theorem} \label{maintheorem}
Under the assumption \refe{regularity-1}-\refe{regularity-2},
there exist two positive constants
$h_0$ and $\tau_0$ such that when $h<h_0$ and $\tau<\tau_0$, 
the FEM systems \refe{fem1}-\refe{fem2} and \refe{fem1-d2}-\refe{fem2-d2} 
are uniquely solvable and the following error estimate holds
\begin{align}
& \max_{0 \leq n \leq N} \Big(
{\| \psi_{h}^{n} - \psi^{n}\|_{L^2}^2}
+{\| \bA_{h}^{n} - \bA^{n}\|_{L^2}^2} \Big)
+ \tau \sum_{m=0}^{N} \| \bm{\sigma}_{h}^{m} - \C \bA^{m}\|_{L^2}^2
\leq C_* (\tau^{2} + h^{2s}) \, ,
\end{align}
where $s = \min \{ r+1, l\}$ and 
$C_*$ is a positive constant independent of $n$, $h$ and $\tau$. 
\end{theorem}

\begin{remark}

The above theorem shows that the convergence rate of the mixed method
depends upon the order of the FEM spaces and also the regularity of 
the exact solution. In \cite{Li_Zhang}, the authors 
proved that on a nonconvex polygon, 
the TDGL equations possess regularity only with 
$\frac{1}{2} < l < \frac{\pi}{\max_{j}\omega_j} < 1$,
where $\max_{j}\omega_j$ denotes the largest interior angle of the polygon. 
More precisely, the two dimensional TDGL equations
admit a solution with $l < \frac{2}{3}$ on an $L$-shape domain.
Therefore, the $L^2$-norm convergence rate of the Galerkin-mixed 
method for the problem 
on an $L$-shape domain is $O(\tau + h^{2/3-\epsilon})$.
We note that the convectional Lagrange FEMs converge to a wrong solution
due to the fact that in space $\bA \in {\bf{H}}^l$ with $l<1$ 
on a nonconvex domain, 
see Example \ref{example2-order} in section \ref{numericalresults}.
Analysis for the TDGL equations in certain three-dimensional geometries 
was given in \cite{Li_Yang}
\end{remark}

In the rest part of this paper, we denote by $C$ a generic positive
constant and $\ee$ a generic small positive constant, 
which are independent of $C_*$, $n$, $h$ and $\tau$. 
We present the Gagliardo--Nirenberg inequality  
in the following lemma which will be 
frequently used in our proofs.
\begin{lemma}
\label{GN} 
{\it ( Gagliardo--Nirenberg inequality \cite{Nirenberg}): 
Let $u$ be a function defined on $\Omega$ in $\mathbb{R}^d$ 
and $\partial ^{s} u$ be any partial derivative of $u$ of order $s$, then
\begin{equation}
\|\partial ^{j} u\|_{L^p} 
\le C \|\partial^{m} u\|_{L^r}^{a} \, \|u\|_{L^q}^{1-a}
+ C \|u\|_{L^q},
\nn
\end{equation}
for $0 \le j < m$ and $\frac{j}{m} \le a \le 1$ with
\[
\frac{1}{p} = 
\frac{j}{d} + a \left( \frac{1}{r} - \frac{m}{d}\right)
+(1-a) \frac{1}{q} \, ,
\]
except $1 < r < \infty$ and $m-j-\frac{n}{r}$ 
is a non-negative integer,
in which case the above estimate 
holds only for $\frac{j}{m} \le a < 1$.
}
\end{lemma}

\section{Preliminaries}
\label{preparation}
\setcounter{equation}{0}

\subsection{An auxiliary elliptic problem}
We consider the elliptic boundary value problem 
\begin{align}
&\C \, \C \, \mathbf{u} - \nabla \, \D \, \mathbf{u} = \mathbf{f} \,
&& \textrm{in $\Omega$},
\label{elliptic-pde} \\
&\C \, \mathbf{u} \times \n = \mathbf{0}, \quad \mathbf{u} \cdot \n = 0
&& \textrm{on $\partial \Omega$},
\label{elliptic-bc}
\end{align}
where $\mathbf{u} = [u_1 \, , u_2 \, , u_3]^{T}$ 
and $\mathbf{f} = [f_1 \, , f_2\, , f_3]^{T}$.
By noting the definition of the standard three dimensional operators
$\C$, $\nabla$ and $\D$,
it is easy to verify that 
${\C \, \C \, \mathbf{u}} - \nabla {\D \, \mathbf{u}} = -\Delta \mathbf{u}$.
By introducing the new variable $\bm{\sigma} = \C \, \mathbf{u}$, 
the mixed formulation of \refe{elliptic-pde}
with the boundary condition \refe{elliptic-bc} is given by
\begin{align}
& \bm{\sigma} = \C \, \mathbf{u} \,
&& \textrm{in $\Omega$},
\label{elliptic-pde-mixed1} \\
&\C \,\bm{\sigma} - \nabla \, \D \, \mathbf{u} = \mathbf{f} \,
&& \textrm{in $\Omega$},
\label{elliptic-pde-mixed2} \\
& \bm{\sigma} \times \n = \mathbf{0} , \quad \mathbf{u} \cdot \n = 0
&& \textrm{on $\partial \Omega$}.
\label{elliptic-bc-mixed}
\end{align}
Then, the weak formulation of the above equation is to find 
$(\bm{\sigma}, \mathbf{u}) \in \accentset{\circ}{\bH}(\C) 
\times \accentset{\circ}{\bH}(\D)$,
such that 
\begin{align}
& (\bm{\sigma} \, , \bm{\chi}) - (\C \, \bm{\chi} \, , \mathbf{u})=0 \, ,
&&\forall \, \bm{\chi} \in  \accentset{\circ}{\bH}(\C)\, ,
\label{elliptic-weak-mixed1} \\
&(\C \,\bm{\sigma} \, , \bv) + (\D \, \mathbf{u} \, , \D \, \bv) 
= (\mathbf{f} \,, \bv) \,,
&&\forall \, \, \bv \in \accentset{\circ}{\bH}(\D) \, .
\label{elliptic-weak-mixed2}
\end{align}
Based on the above mixed weak formulation, the mixed FEM approximation to 
\refe{elliptic-weak-mixed1}-\refe{elliptic-weak-mixed2}
is to find $(\bm{\sigma}_h, \mathbf{u}_h) \in 
\accentset{\circ}{\mathbf{Q}}_h^{r+1} \times \accentset{\circ}{\mathbf{U}}_h^r$,
such that 
\begin{align}
& (\bm{\sigma}_h \, , \bm{\chi}_h) - (\C \, \bm{\chi}_h \, , \mathbf{u}_h)=0 \, ,
&& \forall \, \bm{\chi}_h \in \accentset{\circ}{\mathbf{Q}}_h^{r+1} \, ,
\label{elliptic-fem-mixed1} \\
&(\C \,\bm{\sigma}_{h} \, , \bv_{h}) + (\D \, \mathbf{u}_h \, , \D \, \bv_{h}) 
= (\mathbf{f} \,, \bv_{h}) \,,
&& \forall \, \, \bv_{h} \in  \accentset{\circ}{\mathbf{U}}_h^r \, .
\label{elliptic-fem-mixed2}
\end{align}

Theoretical analysis on convergence 
and stability of the above mixed finite element methods 
can be found in the two seminal papers \cite{AFW, AFW1}, 
while our main concern in this paper is on the 
nonlinear parabolic problem of superconductivity. We summarize the main results
in the following lemma and we refer to \cite{AFW, AFW1} for details.
\begin{lemma}
\label{elliptic-result} 
{\it
Let $(\bm{\sigma}, \mathbf{u})$  and $(\bm{\sigma}_h, \mathbf{u}_h)$ be 
the solution of  
\refe{elliptic-weak-mixed1}-\refe{elliptic-weak-mixed2}
and \refe{elliptic-fem-mixed1}-\refe{elliptic-fem-mixed2},
respectively. Then the following error estimates hold 
\begin{align}
& \|\bm{\sigma} - \bm{\sigma}_h\|_{\mathbf{H}(\C)} 
+ \|\mathbf{u} - \mathbf{u}_h\|_{\mathbf{H}(\D)} 
\nn \\
& \le C \Big(
\inf_{\bm{\chi}_h \in \accentset{\circ}{\mathbf{Q}}_h^{r+1}}
\|\bm{\sigma} - \bm{\chi}_h\|_{\mathbf{H}(\C)}
+\inf_{\mathbf{v}_h \in \accentset{\circ}{\mathbf{U}}_h^r}
\|\mathbf{u} - \mathbf{v}_h\|_{\mathbf{H}(\D)} 
\Big)
\label{elliptic-error1}
\end{align}
}
\end{lemma} 

\subsection{Elliptic and mixed projections}
To prove the optimal error estimates
of the linearized Galerkin-mixed FEM scheme \refe{fem1}-\refe{fem2},
we define a Ritz projection operator 
$R_h : \mathcal{H}^1(\Omega) \rightarrow \mathcal{V}_h^{\widehat{r}}$ 
and a mixed projection operator
${\bf{P}}_h : \mathbf{H}(\C) \times \accentset{\circ}{\mathbf{H}}(\D)
\rightarrow {\mathbf{Q}}_h^{r+1} \times \accentset{\circ}{\mathbf{U}}_h^r$ 
as follows (also see {\cite{thomee}}) : for given $t \in [0,T]$,
find $R_{h}\psi \in \mathcal{V}_{h}^{\widehat{r}}$ and 
${\bf{P}}_h (\bm{\sigma},\bA) := (P_{h}^1 (\bm{\sigma},\bA) \, , P_h^2 (\bm{\sigma},\bA)) 
\in {\mathbf{Q}}_{h}^{r+1} \times \accentset{\circ}{\mathbf{U}}_h^r$
with  $P_{h}^1 (\bm{\sigma},\bA) \times \n =  \pi_h\bm{\sigma} \times \n$ 
on $\partial \Omega$, such that
\begin{eqnarray}
&& \frac{1}{\kappa^2}(\nabla
(\psi-R_{h}\psi) \, , \nabla \omega_{h})
+ M(\psi-R_{h}\psi, \omega_{h})
= 0 \, , \quad \forall \, \omega_{h} \in \mathcal{V}_{h}^{\widehat{r}}
\label{projection-1}
\end{eqnarray}
where $M$ is a positive constant which is chosen to ensure 
the coercivity of \refe{projection-1} and 
\begin{align}
& (\bm{\sigma} - P_{h}^1 (\bm{\sigma},\bA) \, , \bm{\chi}_h) - 
(\C \, \bm{\chi}_h \, , \bA - P_{h}^2 (\bm{\sigma},\bA))=0 \, ,
&& \forall \, \bm{\chi}_h \in \accentset{\circ}{\mathbf{Q}}_h^{r+1} \, ,
\label{projection-1.5} \\
&(\C \,(\bm{\sigma} - P_{h}^1 (\bm{\sigma},\bA)) \, , \bv_{h}) 
+ (\D \, (\bA -P_{h}^2 (\bm{\sigma},\bA))\, , \D \, \bv_{h}) = 0\,,
&& \forall \, \, \bv_{h} \in  \accentset{\circ}{\mathbf{U}}_h^r \, .
\label{projection-2}
\end{align}

We denote the projection error functions by
\begin{equation}
\rho_{\psi} = R_{h} \psi - \psi \, , 
\quad \rho_{\bm{\sigma}} =  P_{h}^1 (\bm{\sigma},\bA) - \bm{\sigma} \, ,
\quad \rho_{\bA} =  P_{h}^2 (\bm{\sigma},\bA) - \bA\, .
\nn
\end{equation}
With the regularity assumption (\ref{regularity-1})-(\ref{regularity-2}),
by standard finite element theory \cite{BS} and Lemma \ref{elliptic-result}, 
we have the following error estimates 
\begin{equation} 
\left \{
\begin{array}{l}
{\| \rho_{\psi} \|_{L^2}} 
+ \| (\rho_{\psi})_{t} \|_{L^2} 
\leq C h^{\min\{2l, \widehat{r}+1\}} ,
\quad {\| \rho_{\psi} \|}_{H^{1}} \leq C h^{\min\{l, \widehat{r}\}} \,,
\\[8pt]
{\| \rho_{\bm{\sigma}} \|_{\mathbf{H}(\C)}} \leq C h^{s}\,, 
\quad
{\| \rho_{\bA} \|_{\mathbf{H}(\D)}} + \| (\rho_{\bA})_{t} \|_{L^2} \leq C h^{s} \,,
\end{array}
\right.
\label{dt1}
\end{equation} 
and the stability result
\begin{eqnarray} 
\| R _h \psi \|_{L^{\infty}}  \le C, \quad
\| R _h \psi \|_{W^{1,3}}  \le C.
\label{dt2} 
\end{eqnarray}

\section{The proof of Theorem {\ref{maintheorem}} }
\label{proofmainresult}
\setcounter{equation}{0}
For $n=0$, $\ldots$, $N$, we denote 
\begin{equation}
e_{\psi}^n = \psi_h^n - {R_h \psi^n} \, , 
\quad e_{\bm{\sigma}}^n = \bm{\sigma}_h^n - P_{h}^1 (\bm{\sigma}^n,\bA^n) \, ,
\quad e_{\bA}^n = \bA_h^n - P_{h}^2 (\bm{\sigma}^n,\bA^n) \, .
\nn
\end{equation}
In this section, we prove that the following inequality holds
for $n = 0$, $\ldots$, $N$
\begin{align}
& {\| e_{\psi}^{n}\|_{L^2}^2}
+ {\| e_{\bA}^{n}\|_{L^2}^2} 
+ \sum_{m=1}^{n} \tau 
\Big( \| e_{\psi}^{m} \|_{H^1}^2
+ \| e_{\bm{\sigma}}^{m}\|_{L^2}^2
+ \| \D \, e_{\bA}^{m}\|_{L^2}^2 \Big)
\leq \frac{C_*}{2} (\tau^{2} + h^{2s})
\label{mainresults}
\end{align}
by mathematical induction. 
Theorem {\ref{maintheorem}} follows immediately from the the projection 
error estimates in \refe{dt1} and the above inequality. 
 
Since 
\begin{align}
& \| e_{\psi}^{0}\|_{L^2}^2 
+ \| e_{\bA}^{0}\|_{L^2}^2 
= \| \pi_h \psi_0 - R_h\psi^0 \|_{L^2}^2
+ \| \pi_h \bA_0 - P_{h}^2 (\bm{\sigma}^0,\bA^0)  \|_{L^2}^2
\le C_1 h^{2s} \, ,
\nn
\end{align}
(\ref{mainresults}) holds for $n = 0$ if we require $ \frac{C_*}{2} \geq C_1$,
we can assume that (\ref{mainresults}) 
holds for $ n \le k-1$ for some $k\ge 1$.
We shall find a constant $C_*$,
which is independent of $n$, $h$, $\tau$,
such that (\ref{mainresults}) holds for $n \leq k$. 
The generic positive constant $C$ in the rest part of this paper 
is independent of $C_*$.

From the mixed variational form \refe{variation1}-\refe{variation2},
the linearized Galerkin-mixed FEM scheme  \refe{fem1}-\refe{fem2},
and the projection \refe{projection-1}-\refe{projection-2}, 
the error functions $e_{\psi}^n$, $e_{\bm{\sigma}}^n$ and $e_{\bA}^n$ satisfy 
\begin{eqnarray}
&&(D_{\tau}e_{\psi}^{n}, \omega_{h}) 
+ \frac{1}{\kappa^2}( \nabla e_{\psi}^{n} \,, \nabla \omega_{h})
+ M ( e_{\psi}^{n} \,, \omega_{h} )
\nn \\
&&  
\quad = -(D_{\tau} \rho_{\psi}^{n}, \omega_{h}) 
+ i\kappa \left(
(\D \,\bA_{h}^{n-1}) \psi_{h}^{n} - (\D \,\bA^{n-1}) \psi^{n}
\,,\omega_{h} \right)
\nn \\
&&  \qquad
+ \frac{i}{\kappa} ({\bA}_{h}^{n-1}{\psi}_{h}^{n} 
- {\bA}^{n-1}{\psi}^{n} \, ,\nabla \omega_{h})
- \frac{i}{\kappa} \big(
(\nabla {\psi}_{h}^{n} \, , {\bA}_{h}^{n-1} \omega_{h})
-(\nabla {\psi}^{n} \, , {\bA}^{n-1} \omega_{h})\big)
\nn \\
&& \qquad
+\left( M \psi_h^n + 
(1-|\psi_{h}^{n-1}|^{2}-|{\bA}_{h}^{n-1}|^2) {\psi}_{h}^{n}
-M \psi^n -
(1-|\psi^{n-1}|^{2}-|{\bA}^{n-1}|^2) {\psi}^{n} \,, \omega_{h} \right)
\nn \\
&& \qquad + R_{\psi}^n 
\nn \\
&& \quad = \sum_{i=1}^{5} J_i^n(\omega_{h}) + R_{\psi}^n \, , 
\quad \forall \, \omega_{h} \in \mathcal{V}_{h}^{\widehat{r}}, 
\quad n=1, 2, \ldots, N
\label{error1}
\end{eqnarray} 
and 
\begin{eqnarray} 
&& (e_{\bm{\sigma}}^n, \bm{\chi}_h) - (\C \, \bm{\chi}_h, e_{\bA}^n) = 0, 
\quad \forall  \, 
\bm{\chi}_h \in \accentset{\circ}{\mathbf{Q}}^{r+1}_{h}\, , 
\quad n=1, 2, \ldots, N,
\label{error1.5}
\\[5pt]
&&(D_{\tau} e_{\bA}^{n},\bv_{h}) 
+ ({\D \, e_{\bA}^{n}} \, ,{\D \, \bv_{h}})
+ ( \C \, e_{\bm{\sigma}}^n \, , \bv_{h})
\nn \\
&& \quad  = -(D_{\tau} \rho_{\bA}^{n},\bv_{h}) 
- \frac{i}{2 \kappa} 
\left[ 
\left( ({\psi}_{h}^{n-1})^{*} 
\nabla {\psi}_{h}^{n-1} 
- {\psi}_{h}^{n-1} \nabla ({\psi}_{h}^{n-1})^{*}, \bv_{h} \right) 
\right.
\nn \\
&& \qquad \left.
-\left( ({\psi}^{n-1})^{*} \nabla {\psi}^{n-1} 
- {\psi}^{n-1} \nabla ({\psi}^{n-1})^{*}, \bv_{h} \right) 
\right]
\nn \\
&& \qquad - (|{\psi}_{h}^{n-1}|^{2} {\bA}_{h}^{n}
- |{\psi}^{n-1}|^{2} {\bA}^{n},\bv_{h}) + R_{\bA}^n
\nn \\
&& \quad = \sum_{i=6}^{8} J_i^n(\bv_{h})  + R_{\bA}^n \, , 
\quad \forall \, \bv_{h} \in \accentset{\circ}{\mathbf{U}}_{h}^{r}, 
\quad n=1, 2, \ldots, N,
\label{error2}
\end{eqnarray}
where
\begin{align}
& R_{\psi}^n = ( D_\tau \psi^n -  
\frac{\partial \psi(\cdot,t_n)}{\partial t}, \omega_{h} )
+ i\kappa \left(\D \,(\bA^{n}-\bA^{n-1}) \psi^{n} \,,\omega_{h} \right)
\nn \\
& \quad \quad + \frac{i}{\kappa}(\nabla \psi^n \,, (\bA^{n-1}-\bA^n)\omega_{h})
- \frac{i}{\kappa}((\bA^{n-1}-\bA^n)\psi^n \,, \nabla \omega_{h})
\nn \\
& \quad \quad
+ \left((|\bA^{n-1}|^2-|\bA^n|^2 +|\psi^{n-1}|^2-|\psi^{n}|^2) \psi^n  \, ,
\omega_{h} \right) 
\nn
\end{align}
and 
\begin{align}
& R_{\bA}^n = ( D_\tau \bA^n -  
\frac{\partial \bA(\cdot,t_n)}{\partial t}, \bv_{h} )
- (|\psi^n|^{2} \bA^n - |\psi^{n-1}|^{2} \bA^{n},\bv_{h})
\nn \\
& \quad \quad
-\frac{i}{2 \kappa}(((\psi^n)^{*} \nabla \psi^n 
- \psi^n \nabla (\psi^n)^{*})-((\psi^{n-1})^{*} \nabla \psi^{n-1}
- \psi^{n-1} \nabla (\psi^{n-1})^{*}), \bv_{h}) 
\nn
\end{align}
define the truncation errors. 

We take $\omega_{h} = e_{\psi}^{n}$ and 
$(\bm{\chi}_h \, , \bv_{h}) = (e_{\bm{\sigma}}^{n} \,, e_{\bA}^{n})$
in \refe{error1}-\refe{error2}, respectively.
It is easy to see that, with the regularity assumptions 
\refe{regularity-1}-\refe{regularity-2},
\begin{align}
| R_{\psi}^n | + | R_{\bA}^n |
\le \ee \|\nabla e_{\psi}^{n}\|_{L^2}^2 
+ C \|e_{\psi}^{n}\|_{L^2}^2 + C \|e_{\bA}^{n}\|_{L^2}^2 
+ \ee^{-1} C \tau^2 \, .
\label{trunc-error}
\end{align}

We now estimate 
$J_i^n(e_{\psi}^{n}) $ for $i=1$, $\ldots$, $5$
in \refe{error1} 
and $J_i^n(e_{\bA}^{n}) $ for $i=6,7,8$
in \refe{error2} term by term.
By noting the projection error estimates \refe{dt1},
we have
\begin{eqnarray}
&& |J_{1}^{n}(e_{\psi}^{n})| \le
C {\| e_{\psi}^{n} \|}^{2} + C h^{2s} \, ,
\nn \\
&& |J_{6}^{n}(e_{\bA}^{n})| \le
C {\| e_{\bA}^{n} \|}^{2} + C h^{2s} \, 
\nn
\end{eqnarray}
and 
\begin{eqnarray}
{\mathcal{R \! \it{e}}}(J_{2}^{n}(e_{\psi}^{n}))
& = & 
{\mathcal{R \! \it{e}}}\left(i \kappa \left( (\D \, \bA_h^{n-1}) e_{\psi}^{n}
+ (\D \, \bA_h^{n-1})R_h\psi^{n}-(\D \, \bA^{n-1})\psi^{n}\, , 
e_{\psi}^{n} \right) \right)
\nn \\
& = & 
{\mathcal{R \! \it{e}}}\left( i \kappa((\D \, \bA_h^{n-1})R_h\psi^{n}-(\D \, \bA^{n-1})\psi^{n}\, , 
e_{\psi}^{n})  \right)
\nn \\
& \le & 
C \big|( (\D \, (e_\bA^{n-1}+\rho_\bA^{n-1})) R_h\psi^{n} 
\, ,e_{\psi}^{n})\big|
+ C \big|((\D \, \bA^{n-1}) \rho_{\psi}^{n} \,, e_{\psi}^{n} ) \big|
\nn \\
& \le & 
C  (\|\D \, e_\bA^{n-1}\|_{L^2}+\|\D \,\rho_\bA^{n-1}\|_{L^2}) 
\|R_h\psi^{n} \|_{L^{\infty}} \|e_{\psi}^{n}\|_{L^2}
\nn \\
& & + C \|\D \, \bA^{n-1}\|_{L^{\infty}} 
\|\rho_{\psi}^{n-1}\|_{L^2} \| e_{\psi}^{n}\|_{L^2}
\nn \\
& \le & 
\ee \| \D \, e_\bA^{n-1} \|_{L^2}
+ \ee^{-1}C (\|e_{\psi}^{n}\|_{L^2}^2 + h^{2s}) \, ,
\nn
\end{eqnarray}
where ${\mathcal{R \! \it{e}}}(J_i)$ denotes the real part of $J_i$ and also,  
we have noted the fact that
\[
{\mathcal{R \! \it{e}}}\big((i \kappa (\D \, \bA_h^{n-1})e_{\psi}^{n}, e_{\psi}^{n})\big) = 0.
\]

For the term $J_3^n$, we see that  
\begin{align} 
|{\mathcal{R \! \it{e}}}\big(J_{3}^{n}(e_{\psi}^{n})\big) |
& \le C |{\mathcal{I\! \it{m}}}\big(( (e_{\bA}^{n-1}+\rho_{\bA}^{n-1})\psi_h^n 
+ \bA^{n-1}(e_{\psi}^{n}+\rho_{\psi}^{n}) \, , \nabla e_{\psi}^{n})\big)|
\nn \\
& \le C |{\mathcal{I\! \it{m}}}\big(((e_{\bA}^{n-1}+\rho_{\bA}^{n-1})e_{\psi}^n \, , 
\nabla e_{\psi}^{n})\big)|
+ C |{\mathcal{I\! \it{m}}}\big(((e_{\bA}^{n-1}+\rho_{\bA}^{n-1})R_h\psi^n \, , 
\nabla e_{\psi}^{n})\big)|
\nn \\
& \quad + C |{\mathcal{I\! \it{m}}}\big((\bA^{n-1}(e_{\psi}^{n}+\rho_{\psi}^{n}) \, , 
\nabla e_{\psi}^{n})\big)|
\nn \\
& \le C |{{\mathcal{I\! \it{m}}}}\big(( e_{\bA}^{n-1} e_{\psi}^{n} \, , \nabla e_{\psi}^{n})\big)|
+ C\| \rho_{\bA}^{n-1} \|_{L^2} \| e_{\psi}^n \|_{L^{\infty}} \| \nabla e_{\psi}^{n} \|_{L^2}  
 + \ee \| \nabla e_{\psi}^{n} \|_{L^2}^2 
\nn \\
& \quad + \ee^{-1}C(\|e_{\psi}^{n} \|_{L^2}^2 + \|e_{\bA}^{n-1} \|_{L^2}^2 + h^{2s})  
\nn \\ 
& \le C |{\mathcal{I\! \it{m}}}\big(( e_{\bA}^{n-1} e_{\psi}^{n} \, , \nabla e_{\psi}^{n})\big)|
+ Ch^{s-\frac{1}{2}} \| e_{\psi}^n \|_{H^1}^2
 + \ee \| e_{\psi}^{n} \|_{H^1}^2
\nn \\
& \quad + \ee^{-1}C(\|e_{\psi}^{n} \|_{L^2}^2 + \|e_{\bA}^{n-1} \|_{L^2}^2 + h^{2s})  
\nn \\ 
& \le  
C |{\mathcal{I\! \it{m}}}\big(( e_{\bA}^{n-1} e_{\psi}^{n} \, , \nabla e_{\psi}^{n})\big)|
+ \ee \|  e_{\psi}^{n} \|_{H^1}^2
+ \ee^{-1}C(\|e_{\psi}^{n} \|_{L^2}^2 + \|e_{\bA}^{n-1} \|_{L^2}^2 + h^{2s})  
\nn
\end{align} 
where we have used the projection error estimate \refe{dt1} and an inverse inequality 
and required $C h^{s-\frac{1}{2}} \le \ee$.  

Similarly, $J_4^n$ is bounded by
\begin{eqnarray}
|{\mathcal{R \! \it{e}}}\big(J_{4}^{n}(e_{\psi}^{n})\big)| 
& \le & 
C | {\mathcal{I\! \it{m}}}\big((\nabla(e_{\psi}^{n}+\rho_{\psi}^{n})\, ,
{\bA}_{h}^{n-1} e_{\psi}^{n} )\big)|
+ C | {\mathcal{I\! \it{m}}}\big((\nabla \psi^n \, , (e_{\bA}^{n-1}+\rho_{\bA}^{n-1})e_{\psi}^{n})\big)|
\nn \\
& \le & C | {\mathcal{I\! \it{m}}}\big((\nabla(e_{\psi}^{n}+\rho_{\psi}^{n})\, ,
(e_{\bA}^{n-1}+\rho_{\bA}^{n-1}) e_{\psi}^{n} )\big)|
+C | {\mathcal{I\! \it{m}}}\big((\nabla(e_{\psi}^{n}+\rho_{\psi}^{n})\, ,
{\bA}^{n-1} e_{\psi}^{n} )\big)|
\nn \\
&  & + \ee \|  e_{\psi}^{n} \|_{H^1}^2
+ \ee^{-1}C(\|e_{\psi}^{n} \|_{L^2}^2 + \|e_{\bA}^{n-1} \|_{L^2}^2 + h^{2s})
\, . 
\nn 
\end{eqnarray} 
By requiring $C h^{s-\frac{1}{2}} \le \ee$ and using inverse inequalities and
integration by parts, 
\begin{eqnarray} 
| {\mathcal{I\! \it{m}}}\big(\nabla\rho_{\psi}^{n}\, ,
(e_{\bA}^{n-1}+\rho_{\bA}^{n-1}) e_{\psi}^{n} )\big)|
&\le& \|\nabla\rho_{\psi}^{n}\|_{L^3} \|e_{\bA}^{n-1}+\rho_{\bA}^{n-1}\|_{L^2} \|e_{\psi}^{n}\|_{L^6}
\nn \\
&\le& C (h^{s} + \|e_{\bA}^{n-1}\|_{L^2}) \|e_{\psi}^{n}\|_{H^1}
\nn \\
&\le &\ee \| e_{\psi}^{n}\|_{H^1}^2 + C (\| e_{\bA}^{n-1} \|_{L^2}^2 + h^{2s}) \, ,
\nn
\end{eqnarray} 
\[
| {\mathcal{I\! \it{m}}}\big((\nabla e_{\psi}^{n}\, ,\rho_{\bA}^{n-1} e_{\psi}^{n} ) \big)|
\le \|\nabla e_{\psi}^{n}\|_{L^2}
\|\rho_{\bA}^{n-1}\|_{L^2} \| e_{\psi}^{n}\|_{L^{\infty}}
\le C h^{s-\frac{1}{2}} \| e_{\psi}^{n}\|_{H^1}^2
\le  \ee  \| e_{\psi}^{n}\|_{H^1}^2 \, ,
\]
\begin{eqnarray} 
| {\mathcal{I \! \it{m}}}\big((\nabla (e_{\psi}^{n}+\rho_{\psi}^{n})\, ,
{\bA}^{n-1} e_{\psi}^{n} ) \big)|
&=& | {\mathcal{I\! \it{m}}}\big((e_{\psi}^{n}+\rho_{\psi}^{n}\, ,
(\D \,{\bA}^{n-1}) e_{\psi}^{n} ) \big)|
\nn \\
&&+ | {\mathcal{I\! \it{m}}}\big((e_{\psi}^{n}+\rho_{\psi}^{n}\, ,
{\bA}^{n-1} \cdot \nabla e_{\psi}^{n} ) \big)|
\nn \\
&\le &\ee \| e_{\psi}^{n} \|_{H^1}^2
+ \ee^{-1} C \| e_{\psi}^{n} \|_{L^2}^2 + \ee^{-1} C h^{2s} \, .  
\nn
\end{eqnarray} 
Then it follows that 
\begin{align} 
|{\mathcal{R \! \it{e}}}\big(J_{4}^{n}(e_{\psi}^{n})\big)| 
\le  C |{\mathcal{I \! \it{m}}}\big(( \nabla e_{\psi}^{n} \, , e_{\bA}^{n-1} e_{\psi}^{n})\big)|
+ \ee \| e_{\psi}^{n} \|_{H^1}^2 
+ \ee^{-1}C(\|e_{\psi}^{n} \|_{L^2}^2 + \|e_{\bA}^{n-1} \|_{L^2}^2 + h^{2s}) \,.
\nn
\end{align} 

For the term $J_{5}$, we have the bound 
\begin{eqnarray}
{\mathcal{R \! \it{e}}}(J_{5}^{n}(e_{\psi}^{n}))
& \le & 
- {\mathcal{R \! \it{e}}}((|\psi_h^{n-1}|^2+|\bA_h^{n-1}|^2)\psi_h^{n} 
- (|\psi^{n-1}|^2+|\bA^{n-1}|^2)\psi^{n} 
\, , e_{\psi}^{n})
\nn \\
&  & + (M+1) \|e_{\psi}^{n}\|_{L^2}^2 + (M+1)(\rho_{\psi}^{n} \, , e_{\psi}^{n})
 \, . 
\nn
\end{eqnarray}
Since
\begin{align} 
&- {\mathcal{R \! \it{e}}}\big((|\bA_h^{n-1}|^2\psi_h^{n}  - |\bA^{n-1}|^2\psi^{n} \, , e_{\psi}^{n})\big)
\nn \\
& \quad = -(|\bA_h^{n-1}|^2e_{\psi}^{n}, \, e_{\psi}^{n})
- {\mathcal{R \! \it{e}}}\big((|\bA_h^{n-1}|^2 R_h\psi^{n}  - |\bA^{n-1}|^2\psi^{n} \, , e_{\psi}^{n})\big)
\nn \\
&  \quad \le -(|\bA_h^{n-1}|^2e_{\psi}^{n}, \, e_{\psi}^{n})
+ |(|\bA^{n-1}|^2 \rho_{\psi}^{n} \, , e_{\psi}^{n})|
+ |((|\bA_h^{n-1}|^2-|\bA^{n-1}|^2) R_h\psi^{n}\, , e_{\psi}^{n})|
\nn \\
&  \quad \le -(|\bA_h^{n-1}|^2e_{\psi}^{n}, \, e_{\psi}^{n})
+ |((\bA_h^{n-1}+\bA^{n-1})\cdot(e_{\bA}^{n-1}+\rho_{\bA}^{n-1}) R_h\psi^{n}\, , e_{\psi}^{n})|
\nn \\
&  \quad  \quad
+ \ee \|e_{\psi}^{n}\|_{H^1}^2 + \ee^{-1} C h^{2s}
\nn \\
&  \quad \le \ee \|e_{\psi}^{n}\|_{H^1}^2 
+ \ee^{-1} C (\|e_{\bA}^{n-1}\|_{L^2}^2 + h^{2s})
\nn
\end{align} 
and 
\begin{eqnarray}
- {\mathcal{R \! \it{e}}}\big((|\psi_h^{n-1}|^2\psi_h^{n}  - |\psi^{n-1}|^2\psi^{n} \, , e_{\psi}^{n})\big)
& \le & \ee \|e_{\psi}^{n}\|_{H^1}^2 
+ \ee^{-1}C (\|e_{\psi}^{n-1}\|_{L^2}^2 + h^{2s}) \,,
\nn
\end{eqnarray}
the real part of $J_{5}^{n}(e_{\psi}^{n})$ can be bounded by
\begin{eqnarray}
{\mathcal{R \! \it{e}}}(J_{5}^{n}(e_{\psi}^{n}))
& \le & 
 \ee \|e_{\psi}^{n}\|_{H^1}^2 + \ee^{-1}C(\|e_{\psi}^{n-1}\|_{L^2}^2
+ \|e_{\bA}^{n-1}\|_{L^2}^2 + h^{2s}) \,.
\nn
\end{eqnarray}
 
For $J_7^n$, we can see that 
\begin{eqnarray}
J_{7}^{n}(e_{\bA}^{n})
& \le & C |{\mathcal{I \! \it{m}}}\big((({\psi}_h^{n-1})^{*} \nabla {\psi}_h^{n-1} 
- ({\psi}^{n-1})^{*} \nabla {\psi}^{n-1} \, , e_{\bA}^{n} )\big)|
\nn \\
& &
+ C |{\mathcal{I \! \it{m}}}\big(({\psi}_h^{n-1} \nabla ({\psi}_h^{n-1})^{*} 
- {\psi}^{n-1} \nabla ({\psi}^{n-1})^{*} \, , e_{\bA}^{n} )\big)| \,.
\label{term7}
\end{eqnarray}
Moreover, the first term in the right hand side of \refe{term7} is bounded by
\begin{align}
& |{\mathcal{I \! \it{m}}}\big((({\psi}_h^{n-1})^{*} \nabla {\psi}_h^{n-1} 
-({\psi}^{n-1})^{*} \nabla {\psi}^{n-1}\, , e_{\bA}^{n} )\big)|
\nn \\
& \quad = |{\mathcal{I \! \it{m}}}\big((({\psi}_h^{n-1})^{*} \nabla (e_{\psi}^{n-1} + \rho_{\psi}^{n-1})
+(e_{\psi}^{n-1} + \rho_{\psi}^{n-1})^{*} \nabla {\psi}^{n-1}\, , e_{\bA}^{n} )\big)|
\nn \\
& \quad \le C |{\mathcal{I \! \it{m}}}\big((({\psi}_h^{n-1})^{*} \nabla (e_{\psi}^{n-1} + \rho_{\psi}^{n-1})
\, , e_{\bA}^{n} )\big)| 
+ \ee\|e_{\psi}^{n-1}\|_{H^1}^2 + \ee^{-1}C(\|e_{\bA}^{n}\|_{L^2}^2 + h^{2s})
\nn \\
& \quad \le  C |{\mathcal{I \! \it{m}}}\big(((e_{\psi}^{n-1})^{*} \nabla (e_{\psi}^{n-1} + \rho_{\psi}^{n-1})
\, , e_{\bA}^{n} )\big)|
+  C |{\mathcal{I \! \it{m}}}\big(((R_h{\psi}^{n-1})^{*} \nabla (e_{\psi}^{n-1} + \rho_{\psi}^{n-1})
\, , e_{\bA}^{n} )\big)|
\nn \\
& \qquad + \ee\|e_{\psi}^{n-1}\|_{H^1}^2 + \ee^{-1}C(\|e_{\bA}^{n}\|_{L^2}^2 + h^{2s})
\nn \\
& \quad \le C |{\mathcal{I \! \it{m}}}\big(( (e_{\psi}^{n-1})^{*} \nabla e_{\psi}^{n-1} \, , 
e_{\bA}^{n} )\big)|
+ C \|e_{\psi}^{n-1}\|_{L^{\infty}}\|\nabla \rho_{\psi}^{n-1}\|_{L^2}\|e_{\bA}^{n}\|_{L^2}
\nn \\
&  \qquad +  C |{\mathcal{I \! \it{m}}}\big(((R_h{\psi}^{n-1})^{*} {\nabla \rho_{\psi}^{n-1}}\, , 
e_{\bA}^{n} )\big)|
+ \ee \| e_{\psi}^{n-1} \|_{H^1}^2
+ \ee^{-1}C(\|e_{\bA}^{n}\|_{L^2}^2 + \|e_{\psi}^{n-1}\|_{L^2}^2 + h^{2s}) \,. 
\nn \\
& \quad \le C |{\mathcal{I \! \it{m}}}\big(( (e_{\psi}^{n-1})^{*} \nabla e_{\psi}^{n-1} \, , 
e_{\bA}^{n} )\big)|
+ C h^{s-\frac{1}{2}} \|e_{\psi}^{n-1}\|_{H^{1}}\|e_{\bA}^{n}\|_{L^2}
\nn \\
&  \qquad + \ee \|\D \,e_{\bA}^{n}\|_{L^2}^2 + \ee \| e_{\psi}^{n-1} \|_{H^1}^2
+ \ee^{-1}C(\|e_{\bA}^{n}\|_{L^2}^2 + \|e_{\psi}^{n-1}\|_{L^2}^2 + h^{2s}) \,. 
\nn \\
& \quad \le C |{\mathcal{I \! \it{m}}}\big(( (e_{\psi}^{n-1})^{*} \nabla e_{\psi}^{n-1} \, , 
e_{\bA}^{n} )\big)|
+ \ee \|\D \,e_{\bA}^{n}\|_{L^2}^2 + \ee \| e_{\psi}^{n-1} \|_{H^1}^2
\nn \\
& \qquad + \ee^{-1}C(\|e_{\bA}^{n}\|_{L^2}^2 + \|e_{\psi}^{n-1}\|_{L^2}^2 + h^{2s}) \,,
\nn 
\end{align}
where we have used an inverse inequality and required that 
$C h^{s-\frac{1}{2}} \le 1$.
Similarly, the second term in the right hand side of \refe{term7} is bounded by 
\begin{align}
& C |{\mathcal{I \! \it{m}}}\big(({\psi}_h^{n-1}\nabla ({\psi}_h^{n-1})^{*}  
-{\psi}^{n-1} \nabla ({\psi}^{n-1})^{*}\, , e_{\bA}^{n} )\big)|
\nn \\
&\le C |{\mathcal{I \! \it{m}}}\big(( e_{\psi}^{n-1} \nabla (e_{\psi}^{n-1})^{*} \, , 
e_{\bA}^{n} )\big)|
+ \ee \|\D \,e_{\bA}^{n}\|_{L^2}^2 + \ee \| e_{\psi}^{n-1} \|_{H^1}^2
\nn \\
& \qquad + \ee^{-1}C(\|e_{\bA}^{n}\|_{L^2}^2 + \|e_{\psi}^{n-1}\|_{L^2}^2 + h^{2s}) \,,
\nn
\end{align}
It follows that
\begin{eqnarray}
J_{7}^{n}(e_{\bA}^{n})
& \le & C |{\mathcal{I \! \it{m}}}\big(( (e_{\psi}^{n-1})^{*} \nabla e_{\psi}^{n-1} \, , 
e_{\bA}^{n} )\big)|
+ \ee \|\D \,e_{\bA}^{n}\|_{L^2}^2 + \ee \| e_{\psi}^{n-1} \|_{H^1}^2
\nn \\
&& + \ee^{-1}C(\|e_{\bA}^{n}\|_{L^2}^2 + \|e_{\psi}^{n-1}\|_{L^2}^2 + h^{2s}) \, .
\nn
\end{eqnarray}

Finally, for the term $J_8^n$, we can derive that
\begin{eqnarray}
J_{8}^{n}(e_{\bA}^{n})
& \le & -(|{\psi}_h^{n-1}|^2\bA_h^{n} - |{\psi}^{n-1}|^2\bA_h^{n} \, , e_{\bA}^{n})
+ C \|e_{\bA}^{n}\|_{L^2}^2 + C h^{2s} 
\nn \\
& = & -(|{\psi}_h^{n-1}|^2 e_{\bA}^{n} \, , e_{\bA}^{n})
-(|{\psi}_h^{n-1}|^2 P_{h}^2 (\bm{\sigma}^n,\bA^n) - |{\psi}^{n-1}|^2\bA_h^{n} \, , e_{\bA}^{n})
+ C \|e_{\bA}^{n}\|_{L^2}^2 + C h^{2s} 
\nn \\
& \le & 
-(|{\psi}_h^{n-1}|^2 e_{\bA}^{n} \, , e_{\bA}^{n})
+ C|(({\psi}_h^{n-1}+{\psi}^{n-1})(e_{\psi}^{n-1}+\rho_{\psi}^{n-1}) 
P_{h}^2 (\bm{\sigma}^n,\bA^n) \, , e_{\bA}^{n})|
\nn \\
& &
+ C \|e_{\bA}^{n}\|_{L^2}^2 + C h^{2s} 
\nn \\
& \le & \ee \|e_{\psi}^{n-1}\|_{H^1}^2 +
\ee^{-1} C (\|e_{\bA}^{n}\|_{L^2}^2 + \|e_{\psi}^{n-1}\|_{L^2}^2 + h^{2s}) \,.
\label{term-8}
\end{eqnarray}

With \refe{trunc-error} and above estimates for $J_i$, 
for $i=1, \ldots, 8$, 
adding \refe{error1.5}, \refe{error2} and
the real part of \refe{error1} together,
we get 
\begin{eqnarray}
&&{\mathcal{R \! \it{e}}}\big((D_{\tau}e_{\psi}^{n}, e_{\psi}^{n}) \big)
+ (D_{\tau} e_{\bA}^{n}, e_{\bA}^{n}) 
+ \frac{1}{\kappa^2} \| \nabla e_{\psi}^{n}\|_{L^2}^2
+ \|\D \, e_{\bA}^{n}\|_{L^2}^2 +  \|e_{\bm{\sigma}}^{n}\|_{L^2}^2 
\nn \\
&& \le C |{\mathcal{I \! \it{m}}}\big(( (e_{\psi}^{n-1})^{*} \nabla e_{\psi}^{n-1} \, , 
e_{\bA}^{n} )\big)|
+C |{\mathcal{I \! \it{m}}}\big(( e_{\bA}^{n-1}e_{\psi}^{n} \, , \nabla e_{\psi}^{n})\big)|
\nn \\
&& \quad
+ \ee ( \|\D \, e_{\bA}^{n}\|_{L^2}^2 + \|\D \, e_{\bA}^{n-1}\|_{L^2}^2
+ \| e_{\psi}^{n}\|_{H^1}^2 + \| e_{\psi}^{n-1}\|_{H^1}^2)
\nn \\
&& \quad
+ \ee^{-1}C (\|e_{\psi}^{n}\|_{L^2}^2
+\|e_{\psi}^{n-1}\|_{L^2}^2 + \|e_{\bA}^{n}\|_{L^2}^2
+ \|e_{\bA}^{n-1}\|_{L^2}^2 + \tau^2 + h^{2s}) \,,
\label{primary-end-1}
\end{eqnarray}
for $n = 1$, $\ldots$, $k$.

In order to estimate 
$|{\mathcal{I \! \it{m}}}\big(( (e_{\psi}^{n-1})^{*} \nabla e_{\psi}^{n-1} \, , e_{\bA}^{n} )\big)|$
and $ |{\mathcal{I \! \it{m}}}\big(( e_{\bA}^{n-1}e_{\psi}^{n} \, , \nabla e_{\psi}^{n})\big)|$
in \refe{primary-end-1},
we use the induction assumption \refe{mainresults} for $ n \le k-1$
and inverse inequalities.

For $\tau \le h^{s}$, we have
\begin{eqnarray}
  |{\mathcal{I \! \it{m}}}\big(( (e_{\psi}^{n-1})^{*} \nabla e_{\psi}^{n-1} \, , e_{\bA}^{n} )\big)|
  & \le & C  \|\nabla e_{\psi}^{n-1}\|_{L^2} \|e_{\psi}^{n-1}\|_{L^6} \|e_{\bA}^{n}\|_{L^3}
  \nn \\
  & \le &  C  \|e_{\psi}^{n-1}\|_{H^1}^2 \,   h^{-\frac{1}{2}} \|e_{\bA}^{n}\|_{L^2}.
  \label{dt-h-prime}
\end{eqnarray}
To estimate $\|e_{\bA}^{n}\|_{L^2}$, from \refe{error2} and estimates \refe{term-8},
we have that
\begin{align}
& {\| e_{\bA}^{n}\|_{L^2}^2} 
+ \tau 
\Big( \| e_{\bm{\sigma}}^{m}\|_{L^2}^2 + \| \D \, e_{\bA}^{m}\|_{L^2}^2 \Big)
\nn \\
&\leq {\| e_{\bA}^{n-1}\|_{L^2}^2}
+ \tau C  \| e_{\psi}^{n-1}\|_{H^1}^2 h^{-\frac{1}{2}}\|e_{\bA}^{n}\|_{L^2}
\nn \\
& \quad
+ \tau(\ee \|e_{\psi}^{n-1}\|_{H^1}^2 +
\ee^{-1} C (\|e_{\bA}^{n}\|_{L^2}^2 + \|e_{\psi}^{n-1}\|_{L^2}^2 + \tau^2+ h^{2s}))
\nn \\
&\leq (\frac{1}{4} + C\tau){\| e_{\bA}^{n}\|_{L^2}^2} + {\| e_{\bA}^{n-1}\|_{L^2}^2}
+ \tau^2 C  \| e_{\psi}^{n-1}\|_{H^1}^4 h^{-1} + C C_* h^{2s}
\nn \\
& \leq \frac{1}{2}{\| e_{\bA}^{n}\|_{L^2}^2} + C C_* h^{2s} \,.
\nn 
\end{align}
where $C\tau \le \frac{1}{4}$.
Then substituting the last estimate into  \refe{dt-h-prime} gives
\begin{eqnarray}
  |{\mathcal{I \! \it{m}}}\big(( (e_{\psi}^{n-1})^{*} \nabla e_{\psi}^{n-1} \, , e_{\bA}^{n} )\big)|
  & \le &  C \sqrt{C_*} h^{s -\frac{1}{2}}  \|e_{\psi}^{n-1}\|_{H^1}^2 \,  
  \nn \\
  & \le &  \epsilon  \|e_{\psi}^{n-1}\|_{H^1}^2 \,  
\end{eqnarray}
where we require $ C \sqrt{C_*} h^{s -\frac{1}{2}} \le \epsilon$.
Similarly,
\begin{eqnarray}
  C|{\mathcal{I \! \it{m}}}\big(( e_{\bA}^{n-1}e_{\psi}^{n} \, , \nabla e_{\psi}^{n})\big)|
  &  \le & C\|e_{\bA}^{n-1}\|_{L^3} \|\nabla e_{\psi}^{n}\|_{L^2} \|e_{\psi}^{n}\|_{L^6}
  \nn \\
  &  \le & C\|e_{\bA}^{n-1}\|_{L^3} \|\nabla e_{\psi}^{n}\|_{L^2}^2
  \nn \\
  &  \le & Ch^{-\frac{1}{2}}\|e_{\bA}^{n-1}\|_{L^2} \|\nabla e_{\psi}^{n}\|_{L^2}^2
  \nn \\
  &  \le & Ch^{-\frac{1}{2}}\sqrt{\frac{C_*}{2} (\tau^{2} + h^{2s})} \|\nabla e_{\psi}^{n}\|_{L^2}^2
  \nn \\
  &  \le & C\sqrt{C_*}h^{s-\frac{1}{2}} \|\nabla e_{\psi}^{n}\|_{L^2}^2
  \nn \\
  &  \le & \epsilon \|\nabla e_{\psi}^{n}\|_{L^2}^2
  \label{dt-h-part1}
\end{eqnarray}
where $h$ satisfies $C\sqrt{C_*}h^{s-\frac{1}{2}} \le \epsilon$.

For $h^{s} \le \tau$, by \refe{mainresults}, we have that for $n \le k-1$
\begin{align}
& {\| e_{\psi}^{n}\|_{L^2}^2} + {\| e_{\bA}^{n}\|_{L^2}^2} 
\leq \frac{C_*}{2} (\tau^{2} + h^{2s}) \leq C_* \tau^{2} \, ,
\nn\\
& \|e_\psi^{n}\|_{H^1}^2 
+\|e_{\bm{\sigma}}^{n}\|_{L^2}^2 + \|\D \, e_\bA^{n}\|_{L^2}^2
\le \tau^{-1} \frac{C_*}{2}(\tau^2 + h^{2s}) \le C_* \tau \, .
\nn
\end{align}
Then, we have
\begin{eqnarray}
  C|{\mathcal{I\! \it{m}}}\big(( (e_{\psi}^{n-1})^{*} \nabla e_{\psi}^{n-1} \, , e_{\bA}^{n} )\big)|
  & \le &  C\|e_{\bA}^{n}\|_{L^3} \| e_{\psi}^{n-1}\|_{H^1} \| e_{\psi}^{n-1}\|_{L^6}
  \nn \\
  & \le & C \sqrt{C_* \tau} (\|e_{\bm{\sigma}}^{n}\|_{L^2}
  + \|\D \, e_\bA^{n}\|_{L^2}) \| e_{\psi}^{n-1}\|_{H^1}
  \nn \\
  & \le & \epsilon (\|e_{\bm{\sigma}}^{n}\|_{L^2}^2 + \|\D \, e_\bA^{n}\|_{L^2}^2
  + \| e_{\psi}^{n-1}\|_{H^1}^2)
  \nn 
\end{eqnarray}
where $\tau$ satisfies $ C C_* \tau \le \epsilon $
and we used the following discrete embedding inequality
(the proof is given in appendix)
\begin{eqnarray}
 \|e_{\bA}^{n}\|_{L^3} 
  & \le &
  C (\|e_{\bm{\sigma}}^{n}\|_{L^2}^2 + \|\D \, e_\bA^{n}\|_{L^2}^2) \,.
  \label{discrete-embedding}
\end{eqnarray}
By noting \refe{dt-h-part1} and \refe{discrete-embedding}, we have
\begin{eqnarray}
  C|{\mathcal{I \! \it{m}}}\big(( e_{\bA}^{n-1}e_{\psi}^{n} \, , \nabla e_{\psi}^{n})\big)|
  & \le & 
  C \|e_{\bA}^{n-1}\|_{L^3} \| e_{\psi}^{n}\|_{H^1}^2
  \nn \\
  & \le &
  C (\|e_{\bm{\sigma}}^{n-1}\|_{L^2} + \|\D \, e_\bA^{n-1}\|_{L^2}) \| e_{\psi}^{n}\|_{H^1}^2
  \nn \\
  & \le &
  C \sqrt{C_*} \tau \| e_{\psi}^{n}\|_{H^1}^2
  \nn \\
  & \le &
  \epsilon  \| e_{\psi}^{n}\|_{H^1}^2
  \nn
\end{eqnarray}
where we require that 
$\tau$ satisfies $ C \sqrt{C_*} \tau \le \epsilon $.

Therefore, from \refe{primary-end-1},
for both cases $\tau \le h^{s}$ and $h^{s} \le \tau$,
we can derive that
\begin{eqnarray}
&&{\mathcal{R \! \it{e}}}\big((D_{\tau}e_{\psi}^{n}, e_{\psi}^{n}) \big)
+ (D_{\tau} e_{\bA}^{n}, e_{\bA}^{n}) 
+ \frac{1}{\kappa^2} \| \nabla e_{\psi}^{n}\|_{L^2}^2
+ M \|e_{\psi}^{n}\|_{L^2}^2 
+ \|\D \, e_{\bA}^{n}\|_{L^2}^2 +  \|e_{\bm{\sigma}}^{n}\|_{L^2}^2 
\nn \\
&& \quad \le  \ee ( \| e_{\bm{\sigma}}^{n}\|_{L^2}^2 + \|\D \, e_{\bA}^{n}\|_{L^2}^2
+ \| e_{\psi}^{n}\|_{H^1}^2 + \| e_{\psi}^{n-1}\|_{H^1}^2)
\nn \\
&& \qquad
+ \ee^{-1}C (\|e_{\psi}^{n}\|_{L^2}^2
+\|e_{\psi}^{n-1}\|_{L^2}^2 + \|e_{\bA}^{n}\|_{L^2}^2
+ \|e_{\bA}^{n-1}\|_{L^2}^2 + \tau^2 + h^{2s}) \, .
\nn
\end{eqnarray}
By choosing a small $\ee$ and summing up the last inequality
from $1$ to $n$, we arrive at
\begin{align} 
\| e_{\psi}^{n} \|_{L^2}^2 & +  \| e_{\bA}^{n} \|_{L^2}^2 
+ \tau \sum_{m=1}^{n} \left( \| e_{\psi}^{m} \|_{H^1}^2
+ \|\D \, e_{\bA}^{m}\|_{L^2}^2 +  \|e_{\bm{\sigma}}^{m}\|_{L^2}^2 \right)
\nn \\
& \le
C \tau \sum_{m=1}^{n} \left( \|e_{\psi}^{m}\|_{L^2}^2 + \|e_{\bA}^{m}\|_{L^2}^2 \right) 
+ C ( \tau^2 + h^{2s}) \, .
\nn
\end{align} 
With the help of the Gronwall's inequality,
we can deduce that 
\begin{align}
& {\| e_{\psi}^{n}\|_{L^2}^2} + {\| e_{\bA}^{n}\|_{L^2}^2} 
+ \sum_{m=1}^{n} \tau \Big( \| e_{\psi}^{m} \|_{H^1}^2 + \| e_{\bm{\sigma}}^{m}\|_{L^2}^2
+ \| \D \, e_{\bA}^{m}\|_{L^2}^2 \Big)
\leq C (\tau^{2} + h^{2s}) \,.
\nn
\end{align}
Thus \refe{mainresults} holds for $n=k$, 
if we take $C_*>2C$.
The induction is closed and the proof of Theorem \ref{maintheorem} is complete. 
\quad \endproof

\section{Numerical results}
\label{numericalresults}
\setcounter{equation}{0}

In this section, we provide several numerical 
examples in both two and three dimensional spaces
to confirm our theoretical  analysis and show the efficiency of the 
linearized Galerkin-mixed FEM scheme. 
The computations are carried out with the free software FEniCS {\cite{fenics}}.

\subsection{Two dimensional numerical experiments}

\vspace{0.1 in}
\begin{example}
\rm
\label{example1}
In this example, we consider the following two dimensional artificial problem
\begin{align} 
& \frac{\partial \psi}{\partial t} 
- i \kappa ({\D \, \bA}) \psi 
+ (\frac{i}{\kappa} \nabla + \bA)^{2} \psi + (|\psi|^{2}-1) \psi = g \, ,
\label{pdetest1}\\
&\frac{\partial \bA}{\partial t} 
- \nabla {\D \, \bA} + {\C \, \curl \, \bA} 
+ \frac{i}{2 \kappa}(\psi^{*} \nabla \psi - \psi \nabla \psi^{*}) 
+ |\psi|^{2} \bA = \C \, H_e + \mathbf{f} \, , 
\label{pdetest2}
\end{align} 
for $t \in (0,T]$, $x \in \Omega$, with boundary and initial conditions
\begin{align} 
& \frac{\partial \psi}{\partial n} = 0, 
\quad \curl \, \bA = H_e, \quad \bA \cdot \mathbf{n} = 0,
&& \mathrm{on}\ \partial \Omega ,
\nn
\\ 
&\psi(x,0) = \psi_{0}(x), \quad \bA(x,0) = \bA_{0}(x), 
&& \mathrm{in}\ \Omega ,
\nn 
\end{align} 
where $\Omega = (0,1) \times (0,1)$ 
and $ \kappa = 1$. 
The functions $\mathbf{f}$, $g$, $\psi_0$ and $\bA_0$ are chosen 
correspondingly to the exact solution 
\begin{align}
\psi = \exp{(-t)} (\cos(\pi x) + i \cos(\pi y)) \, , 
\quad \bA = \left[
\begin{array}{c}
\exp(y-t)\sin(\pi x) \\
\exp(x-t)\sin(\pi y)
\end{array}
\right ]
\nn
\end{align}
with 
\begin{equation}
H_e = \exp(x-t)\sin(\pi y) - \exp(y-t)\sin(\pi x) \, .
\nn
\end{equation}
We set $T = 1.0$ in this example.

\begin{figure}[htp]
\vspace{5pt}
\centering
\begin{tabular}{c}
\epsfig{file=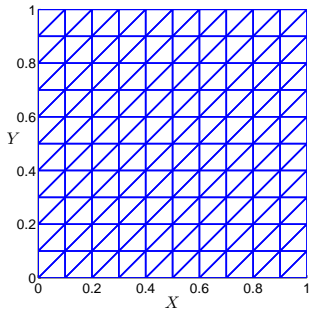,width=2.3in}
\end{tabular}
\caption{A uniform mesh on the unit square with $M=10$ (Example \ref{example1}).}
\label{mesh}
\end{figure}

We use a uniform triangular partition with $M+1$ vertices 
in each direction, see Figure \ref{mesh} for $M=10$, 
where $h = \frac{\sqrt{2}}{M}$.
We solve the system \refe{pdetest1}-\refe{pdetest2} 
by the proposed linearized backward
Euler Galerkin-mixed finite element methods \refe{fem1-d2}-\refe{fem2-d2} 
with $(\mathcal{V}_h^{\widehat{r}}, \mathbf{Q}_h^{r+1} \times \mathbf{U}_h^{r} )$
for $r=0$, $1$, $2$.
As $\sigma = \curl \, \bA$ is a real scalar function in the two dimensional case,
correspondingly its approximation space $\mathbf{Q}_h^{r+1}$ 
is the Lagrange element space consisting of all
piecewise polynomials up to $(r+1)$-th order.
As the exact solution is smooth,to confirm our $L^{2}$-norm error estimates, 
we set $\tau = \left( \frac{1}{M} \right)^{r+1}$ in our simulation and
we present in Table \ref{rt-d2} the $L^2$-norm errors of $\psi$, $\bA$, $\sigma$.
We can see clearly that the $L^{2}$-norm errors  
are proportional to $h^{r+1}$, $r=0$, $1$, $2$, 
which confirm the optimal convergence of the scheme in two dimensional space.
{\setlength{\extrarowheight}{3pt}
\begin{table}[htp]
\begin{center}
\caption{$L^{2}$-norm errors of $\psi$, $\bA$
and $\sigma$ on the unit square (Example \ref{example1})..}
\label{rt-d2}
\begin{tabular}{c|ccc}
\hline
\hline
$(\mathcal{V}_h^{1}, \mathbf{Q}_h^{1} \times \mathbf{U}_h^{0} )$ 
~~ $\tau = \frac{1}{M}$
    & $\|\psi_{h}^N-\psi(\cdot,t_N)\|_{L^2}$  &
      $\|\bA_{h}^N-\bA(\cdot,t_N)\|_{L^2}$    &
      $\|\sigma_{h}^N-\sigma(\cdot,t_N)\|_{L^2}$   \\
\hline
$M=64 $ & 3.1216e-02  & 1.0296e-02  & 2.0804e-03 \\
$M=128$ & 1.6284e-02  & 5.1766e-03  & 1.0637e-03 \\
$M=256$ & 8.3156e-03  & 2.5959e-03  & 5.3770e-04 \\
\hline
order   & 0.95        & 0.99        & 0.98       \\
\hline
\hline
$(\mathcal{V}_h^{1}, \mathbf{Q}_h^{2} \times \mathbf{U}_h^{1} )$ 
~~ $\tau = \frac{1}{M^2}$
    & $\|\psi_{h}^N-\psi(\cdot,t_N)\|_{L^2}$  &
      $\|\bA_{h}^N-\bA(\cdot,t_N)\|_{L^2}$    &
      $\|\sigma_{h}^N-\sigma(\cdot,t_N)\|_{L^2}$   \\
\hline
$M=16$ & 9.9391e-03  & 8.1354e-04  & 6.3113e-04 \\
$M=32$ & 2.4956e-03  & 2.0431e-04  & 1.5807e-04 \\
$M=64$ & 6.2454e-04  & 5.1132e-05  & 3.9533e-05 \\
\hline
order  & 2.00        & 2.00        & 2.00       \\
\hline
\hline
$(\mathcal{V}_h^{2}, \mathbf{Q}_h^{3} \times \mathbf{U}_h^{2} )$ 
~~ $\tau = \frac{1}{M^3}$
    & $\|\psi_{h}^N-\psi(\cdot,t_N)\|_{L^2}$  &
      $\|\bA_{h}^N-\bA(\cdot,t_N)\|_{L^2}$    &
      $\|\sigma_{h}^N-\sigma(\cdot,t_N)\|_{L^2}$   \\
\hline
$M=8$ & 4.1680e-03   & 4.5426e-04  & 2.7586e-04 \\
$M=16$& 5.2687e-04   & 5.7080e-05  & 3.4243e-05 \\
$M=32$& 6.6130e-05   & 7.1400e-06  & 4.2634e-06 \\
\hline
order & 2.99         & 3.00        & 3.01       \\
\hline
\hline
\end{tabular}
\end{center}
\end{table}
}

To show the unconditional convergence of the methods, 
we solve \refe{pdetest1}-\refe{pdetest2} by the scheme \refe{fem1-d2}-\refe{fem2-d2} 
with $(\mathcal{V}_h^{1}, \mathbf{Q}_h^{2} \times \mathbf{U}_h^{1} )$ and  
three different time steps $\tau=0.1$, $0.01$, $0.001$. 
For each fixed $\tau$, we take $M=8$, $16$, $32$, $64$, $128$. 
The $L^2$ errors of $\psi$, $\bA$
and $\sigma$ are presented in Figure \ref{stab_RT}.
From Figure \ref{stab_RT}, we can see that for each fixed $\tau$, when  
the mesh is refined gradually, each $L^2$ error converges to a small constant 
of $O(\tau)$, which shows clearly 
that the proposed scheme \refe{fem1-d2}-\refe{fem2-d2} 
is unconditionally stable and no time step restriction is needed.
\begin{figure}[htp]
\vspace{5pt}
\centering
\begin{tabular}{ccc}
\epsfig{file= 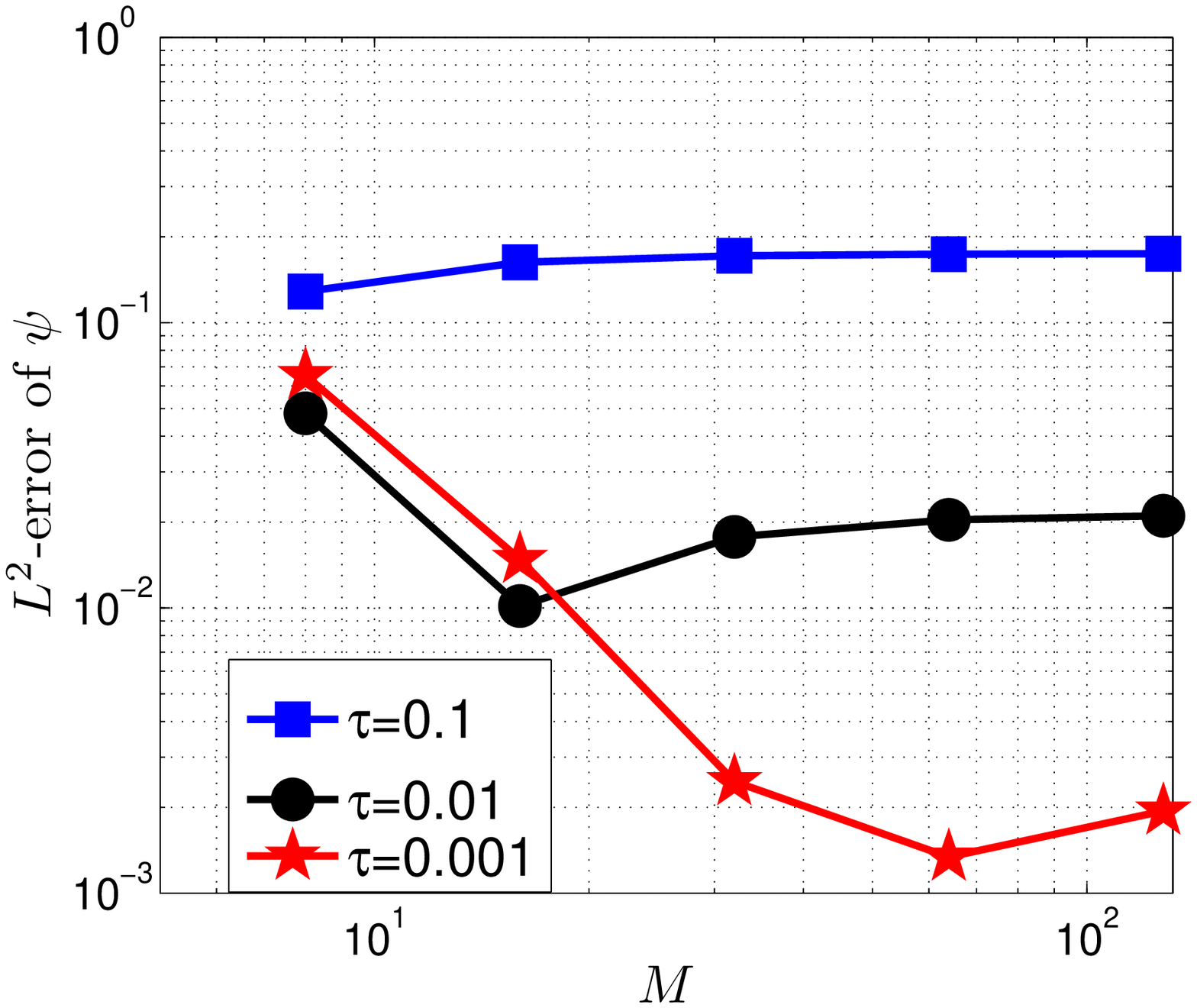,  width=1.9in,height=1.8in} &
\epsfig{file= 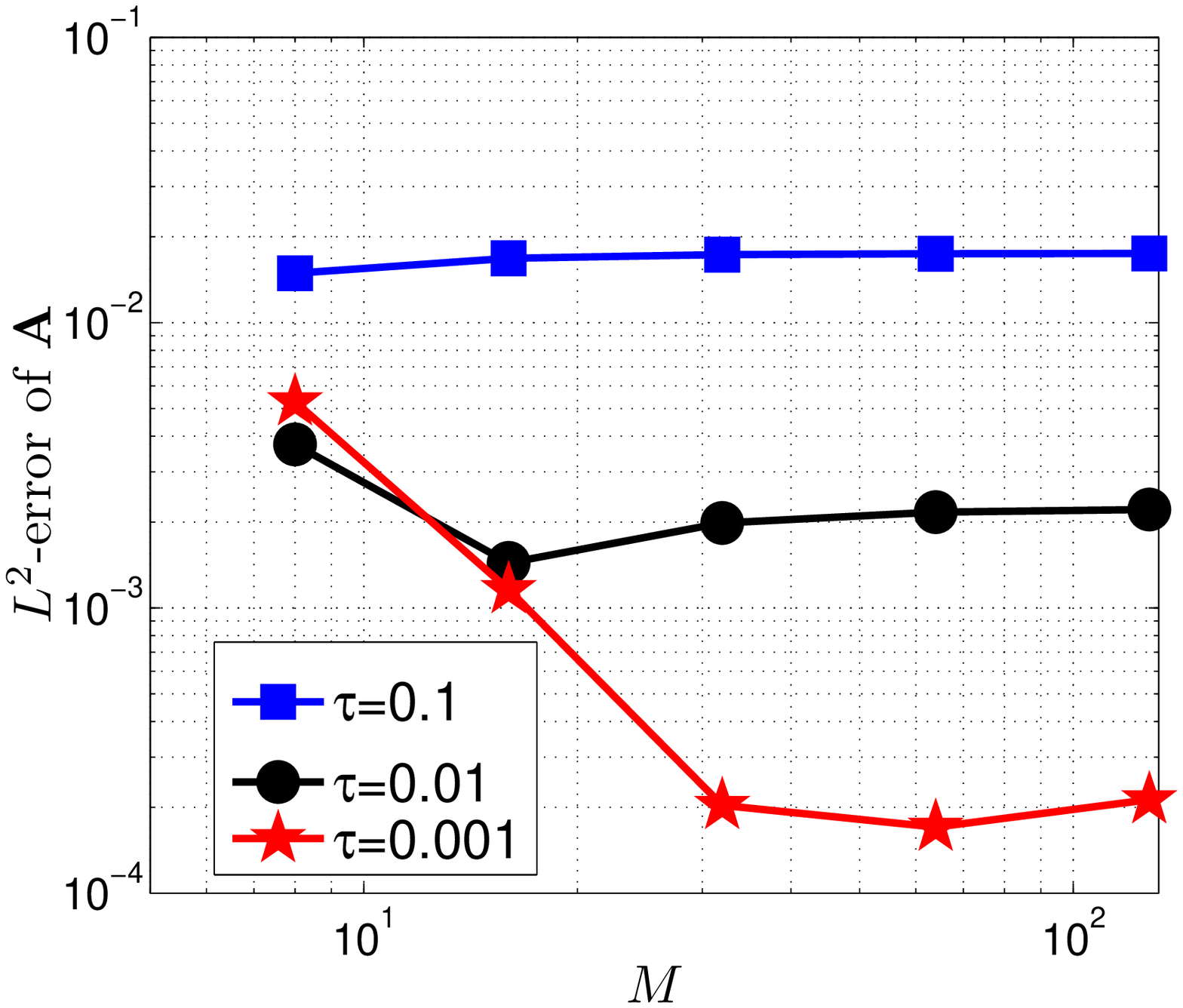,  width=1.9in,height=1.8in} &
\epsfig{file= 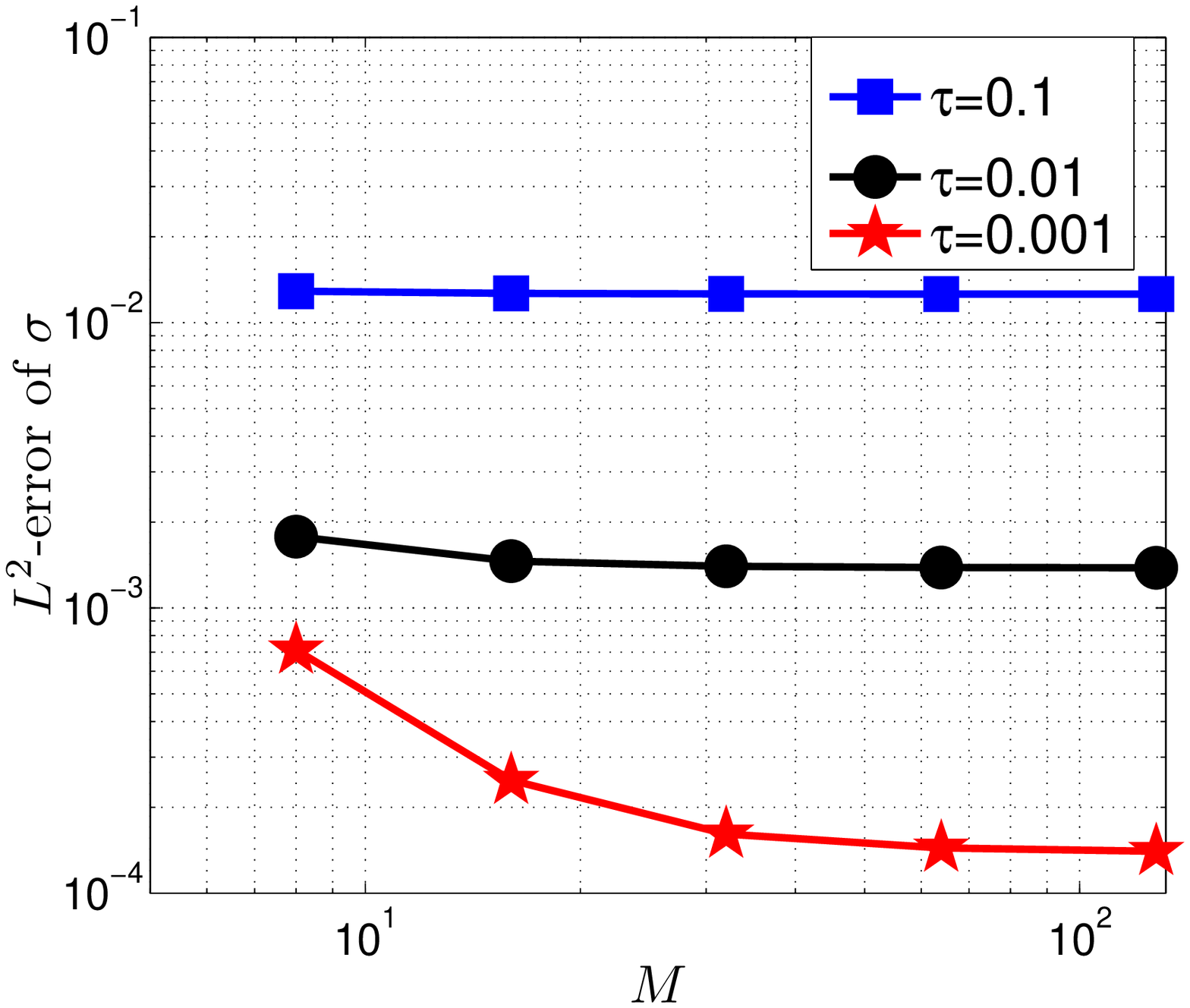,width=1.9in,height=1.8in} 
\end{tabular}
\caption{$L^2$ errors of $\psi$, $\bA$ and $\sigma$ with 
$(\mathcal{V}_h^{1}, \mathbf{Q}_h^{2} \times \mathbf{U}_h^{1} )$ 
(Example \ref{example1}).}
\label{stab_RT}
\end{figure}
\end{example}

\vspace{0.1 in}
\begin{example}
\rm
\label{example2-order}
In this example, we investigate the TDGL equations 
\refe{pdetest1}-\refe{pdetest2} on an $L$-shape domain 
with a nonsmooth solution by the Galerkin-mixed method. 
This example was studied in \cite{Li_Zhang} by a projected FEM based 
on Hodge decomposition.
Here, we use the same exact solution as in \cite{Li_Zhang}, i.e.,
$ \kappa = 10$, and the functions $\mathbf{f}$, $g$, 
$H_e= \curl \, \bA$, $\psi_0$ and $\bA_0$ are chosen 
correspondingly to the following $\psi$ and $\bA$
\begin{align}
 \psi = t^2 \Phi(r) r^{\frac{2}{3}} \cos \left(\frac{2 \theta}{3} \right),
\quad \bA = \left[
\begin{array}{c}
\left( \frac{4}{3} t^2 \Phi(r) r^{-\frac{1}{3}} + t^2 \Phi'(r) r^{\frac{2}{3}} \right)
\cos \left(\frac{\theta}{3} \right)
 \\[4pt]
\left( \frac{4}{3} t^2 \Phi(r) r^{-\frac{1}{3}} + t^2 \Phi'(r) r^{\frac{2}{3}} \right)
\sin \left(\frac{\theta}{3} \right)
\end{array}
\right ] \, ,
\nn
\end{align}
where $(r,\theta)$ denotes the two-dimensional polar coordinates.
$\Phi(s)$ in the above expressions is a cut-off function defined by
\begin{equation}
\Phi(s) = \left \{
\begin{array}{ll}
0.1 & \textrm{if $s < 0.1$},
\\[4pt]
\Upsilon(s) & \textrm{if $0.1 \le s \le 0.4$},
\\[4pt]
0 & \textrm{if $s > 0.4$},
\end{array}
\right.
\nn
\end{equation}
where $\Upsilon(s)$ is a seventh order polynomial satisfying 
\begin{equation}
\left \{
\begin{array}{l}
\Upsilon(0.1)=0.1,\quad \Upsilon(0.4)=0,
\\[4pt]
\Upsilon^{(1)}(0.1)=\Upsilon^{(2)}(0.1)=\Upsilon^{(3)}(0.1)=0,
\\[4pt]
\Upsilon^{(1)}(0.4)=\Upsilon^{(2)}(0.4)=\Upsilon^{(3)}(0.4)=0.
\end{array}
\right.
\nn
\end{equation}
It is noted that in spatial space, 
the above exact solution only has the regularity
\[
\psi \in \mathcal{H}^{\frac{5}{3}-\epsilon },\quad
\bA \in \mathbf{H}^{\frac{2}{3}-\epsilon }, \quad
\sigma \in {H}^{\frac{5}{3}-\epsilon },\quad
\D \, \bA \in {H}^{\frac{5}{3}-\epsilon },\quad 
\quad \textrm{for $\epsilon > 0$}\, .
\]

\begin{figure}[htp]
\vspace{5pt}
\centering
\begin{tabular}{c}
\epsfig{file=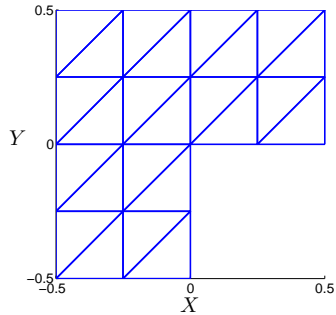,width=2.3in}
\end{tabular}
\caption{A uniform mesh on the $L$-shape domain with $M=4$
(Example \ref{example2-order}).}
\label{mesh-lshape}
\end{figure}

Now, we use the Galerkin-mixed method \refe{fem1-d2}-\refe{fem2-d2} to
solve the above problem on uniform meshes as shown in Figure \ref{mesh-lshape},
where $h=\frac{\sqrt{2}}{M}$. 
We set the the terminal time $T = 1.0$ and do computation 
with $(\psi_h, \bA_h,  \sigma_h) \in 
(\mathcal{V}_h^{1}, \mathbf{Q}_h^{1}, \mathbf{U}_h^{0} )$
and the time step $\tau = \frac{1}{M}$. 
The $L^2$ errors are shown in Table \ref{order-l-shape-example2}.
From Table \ref{order-l-shape-example2}, 
we can see that the convergence rates of 
the Galerkin-mixed method for the three components $\psi$, $\bA$ and $\sigma$,
are in the order of $O(h^{1.14})$, $O(h^{0.86})$ and $O(h^{1.91})$, 
respectively, which indicates that our method is suitable 
for problems on nonconvex domains.

{\setlength{\extrarowheight}{3pt}
\begin{table}[htp]
\begin{center}
\caption{$L^{2}$ norm errors of $\psi$, $\bA$
and ${\sigma}$ by $(\mathcal{V}_h^{1}, 
\mathbf{Q}_h^{1} \times \mathbf{U}_h^{0} )$
on the $L$-shape domain (Example \ref{example2-order}).}
\label{order-l-shape-example2}
\begin{tabular}{c|ccc}
\hline
\hline
$(\mathcal{V}_h^{1}, \mathbf{Q}_h^{1} \times \mathbf{U}_h^{0})$
~~ $\tau = \frac{1}{M}$
    & $\|\psi_{h}^N-\psi(\cdot,t_N)\|_{L^2}$  &
      $\|\bA_{h}^N-\bA(\cdot,t_N)\|_{L^2}$    &
      $\|{\sigma}_{h}^N-{\sigma}(\cdot,t_N)\|_{L^2}$   \\
\hline
$M= 32$ & 6.5548e-04 & 1.7728e-02 & 1.2104e-01 \\
$M= 64$ & 2.6140e-04 & 9.2365e-03 & 3.3953e-02 \\
$M=128$ & 1.2130e-04 & 5.0878e-03 & 8.9061e-03 \\
$M=256$ & 5.9174e-05 & 2.9632e-03 & 2.2696e-03 \\
\hline     
order   & 1.04       & 0.86       & 1.91       \\
\hline
\hline
\end{tabular}
\end{center}
\end{table}
}

For comparison, we also solve this artificial problem 
by a Lagrange type FEM under the same data settings. 
The linearized backward Euler Lagrange FEM
is to find $\psi_{h}^{n} \in \mathcal{V}_{h}^{1}$ and 
$\bA_{h}^{n} \in \accentset{\circ}{\mathbf{V}}_h^{1}$,
such that for $n = 1,2, \ldots, N$, 
\begin{eqnarray}
&&(D_{\tau}\psi_{h}^{n}, \omega_{h})
- i\kappa ((\D \,\bA_{h}^{n-1}) \psi_{h}^{n},\omega_{h}) 
+ \left((\frac{i}{\kappa}\nabla + \bA_h^{n-1}) {\psi}_{h}^{n} \,, 
(\frac{i}{\kappa}\nabla + \bA_h^{n-1})  \omega_{h} \right)
\nn \\
&&
\qquad
- ((| \bA_h^{n-1} |^{2}+ |\psi_{h}^{n-1}|^{2}-1)\psi_{h}^{n}, \omega_{h})
= 0
\quad \forall \, \omega_{h} \in \mathcal{V}_{h}^{1} \, , 
\label{fem1-d2-old}
\end{eqnarray}
and
\begin{eqnarray}
&&(D_{\tau}\bA_{h}^{n},\bv_{h}) 
+ ({\D \, {\bA}_{h}^{n}} \, ,{\D \, \bv_{h}})
+ ({\curl \, {\bA}_{h}^{n}} \, ,{\curl \, \bv_{h}})
+ (|{\psi}_{h}^{n-1}|^{2} {\bA}_{h}^{n},\bv_{h}) 
\nn \\
&& \qquad  = ( {H}_e^{n} \, , \curl \, \bv_{h}) 
 - \frac{i}{2 \kappa} \left( ({\psi}_{h}^{n-1})^{*} 
\nabla {\psi}_{h}^{n-1} 
- {\psi}_{h}^{n-1} \nabla ({\psi}_{h}^{n-1})^{*}, \bv_{h} \right) 
\quad \forall \, \bv_{h} \in \accentset{\circ}{\mathbf{V}}_{h}^{1} \, ,
\label{fem2-d2-old}
\end{eqnarray}
where $\accentset{\circ}{\mathbf{V}}_{h}^{1}$ denotes the linear Lagrange
FE subspace of 
$\left\{ \bA \, \big| \, \bA \in {\bH}^1, \,
\bA \cdot \mathbf{n} \big|_{\partial \Omega} = 0 \right\} $.
$\psi_{h}^{0} = \pi_h\psi_{0}$ and $\bA_{h}^{0} = \pi_h \bA_0$ are used
at the initial time step.
The $L^2$ errors obtained by the conventional Lagrange FEM 
\refe{fem1-d2-old}-\refe{fem2-d2-old} with $M=32$, $64$, $128$, $256$ 
are presented in Table \ref{fem-order-l-shape-example2}.
We see clearly that the numerical solution of 
the Lagrange FEM does not converge. 
We refer to \cite{GS,Li_Zhang2} for the numerical observation of 
the spurious convergence phenomenon 
by conventional Lagrange FEMs in the vortex motion simulations on an $L$-shape domain.

{\setlength{\extrarowheight}{3pt}
\begin{table}[htp]
\begin{center}
\caption{$L^{2}$ norm errors of $\psi$, $\bA$
and ${\sigma}$ on the $L$-shape domain by the conventional Lagrange FEM
\refe{fem1-d2-old}-\refe{fem2-d2-old}
(Example \ref{example2-order}).}
\label{fem-order-l-shape-example2}
\begin{tabular}{c|ccc}
\hline
\hline
$(\mathcal{V}_h^{1},  \accentset{\circ}{\mathbf{V}}_{h}^{1})$
~~ $\tau = \frac{1}{M}$
    & $\|\psi_{h}^N-\psi(\cdot,t_N)\|_{L^2}$  &
      $\|\bA_{h}^N-\bA(\cdot,t_N)\|_{L^2}$    &
      $\|{\sigma}_{h}^N-{\sigma}(\cdot,t_N)\|_{L^2}$   \\
\hline
$M= 32$ & 3.8735e-03 & 9.6105e-02 & 5.1970e-01 \\
$M= 64$ & 3.8616e-03 & 9.8059e-02 & 4.7659e-01 \\
$M=128$ & 3.8466e-03 & 1.0298e-01 & 4.6553e-01 \\
$M=256$ & 3.8395e-03 & 1.0539e-01 & 4.6204e-01 \\
\hline
\hline
\end{tabular}
\end{center}
\end{table}
}

\end{example}

\subsection{Three dimensional numerical experiments}
Numerical results for a three-dimensional artificial example was reported 
in \cite{GS}. For completeness, we include some similar results 
in this section to confirm the convergence rate of the Galerkin-mixed 
method for the problems in three-dimensional space.

\vspace{0.1 in}
\begin{example}
\rm
\label{example3d}
We consider the following three dimensional artificial problem
\begin{align} 
& \frac{\partial \psi}{\partial t} 
- i  \kappa ({\D \, \bA}) \psi 
+ (\frac{i}{\kappa} \nabla + \bA)^{2} \psi + (|\psi|^{2}-1) \psi = g \, ,
\label{pdetest1-3d}\\
&\frac{\partial \bA}{\partial t} 
- \nabla {\D \, \bA} + {\C \, \C \, \bA} 
+ \frac{i}{2 \kappa}(\psi^{*} \nabla \psi - \psi \nabla \psi^{*}) 
+ |\psi|^{2} \bA = \C \, \mathbf{H}_e + \mathbf{f} \, , 
\label{pdetest2-3d}
\end{align} 
for $t \in (0,T]$, $x \in \Omega$, with 
\begin{align} 
& \frac{\partial \psi}{\partial \mathbf{n}} = 0, 
\quad \C \, \bA \times \mathbf{n} = \mathbf{H}_e \times \mathbf{n}, 
\quad \bA \cdot \mathbf{n} = 0,
&& \mathrm{on}\ \partial \Omega ,
\nn
\\ 
&\psi(x,0) = \psi_{0}(x), \quad \bA(x,0) = \bA_{0}(x), \quad 
&& \mathrm{in}\ \Omega ,
\nn 
\end{align} 
where we set $\Omega = (0,1) \times (0,1) \times (0,1)$ and $\kappa = 1$. 
The functions $\mathbf{f}$, $g$, $\psi_0$ and $\bA_0$ are chosen 
correspondingly to the exact solution 
\begin{align}
 \psi = \exp{(t)} \left(
\cos(\pi x)\cos(\pi z) + i  \cos(\pi y) \cos(\pi z) \right) \, , 
\quad \bA = \left[
\begin{array}{c}
\exp(t)\sin(2 \pi x)\sin(2 \pi y) \\
\exp(t)\sin(2 \pi y)\sin(2 \pi z) \\
\exp(t)\sin(2 \pi z)
\end{array}
\right]
\nn
\end{align}
with 
\begin{equation}
\mathbf{H}_e =
 \left[
\begin{array}{c}
-2 \pi \exp(t)\sin(2 \pi y)\cos(2 \pi z) \\
0 \\
-2 \pi \exp(t)\sin(2 \pi x)\cos(2 \pi y) 
\end{array}
\right] \,.
\nn
\end{equation}
We also set the terminate time $T = 1.0$ in this example.

A uniform tetrahedral mesh with $M+1$ vertices in each direction of the cube 
is used in our computation, where $h = \frac{\sqrt{3}}{M}$.
We solve the system by the  Galerkin-mixed FEM scheme \refe{fem1}-\refe{fem2} 
with 
$(\mathcal{V}_h^{\widehat{r}}, \mathbf{Q}_h^{r+1} \times \mathbf{U}_h^{r})$
for $r=0$, $1$, respectively.
Here $\bm{\sigma} = \C \, \bA$ is a 
three dimensional vector function and correspondingly
$\mathbf{Q}_h^{r+1}$ is the $(r+1)$-th order
first type {N\'ed\'elec} element space.
To confirm our error estimates in the $L^{2}$ norm, 
we choose $\tau = \left( \frac{1}{M} \right)^{r+1}$ for 
$(\mathcal{V}_h^{\widehat{r}}, \mathbf{Q}_h^{r+1} \times \mathbf{U}_h^{r})$.
We present in Table \ref{order-3d} the $L^2$-norm errors  
of $\psi$, $\bA$, $\bm{\sigma}$.
We can see clearly that the $L^{2}$ norm errors  
are in the order of $O(h^{r+1})$, $r=0$, $1$, 
which confirm our error analysis in three dimensional space.
{\setlength{\extrarowheight}{3pt}
\begin{table}[htp]
\begin{center}
\caption{$L^{2}$ norm errors of $\psi$, $\bA$
and $\bm{\sigma}$ on the unit cube (Example \ref{example3d}).}
\label{order-3d}
\begin{tabular}{c|ccc}
\hline
\hline
$(\mathcal{V}_h^{1}, \mathbf{Q}_h^{1} \times \mathbf{U}_h^{0})$
~~ $\tau = \frac{1}{M}$
    & $\|\psi_{h}^N-\psi(\cdot,t_N)\|_{L^2}$  &
      $\|\bA_{h}^N-\bA(\cdot,t_N)\|_{L^2}$    &
      $\|\bm{\sigma}_{h}^N-\bm{\sigma}(\cdot,t_N)\|_{L^2}$   \\
\hline
$M=8 $ & 2.0951e-01  & 7.0869e-01  & 3.9600e+00   \\
$M=16$ & 8.3881e-02  & 3.3623e-01  & 1.9502e+00   \\
$M=32$ & 3.7858e-02  & 1.6561e-01  & 9.7032e-01   \\
\hline   
order  & 1.234  & 1.049  & 1.015   \\
\hline
\hline
$(\mathcal{V}_h^{1}, \mathbf{Q}_h^{2} \times \mathbf{U}_h^{1})$
~~ $\tau = \frac{1}{M^2}$
    & $\|\psi_{h}^N-\psi(\cdot,t_N)\|_{L^2}$  &
      $\|\bA_{h}^N-\bA(\cdot,t_N)\|_{L^2}$    &
      $\|\bm{\sigma}_{h}^N-\bm{\sigma}(\cdot,t_N)\|_{L^2}$   \\
\hline
$M=4 $ & 4.2803e-01  & 3.3717e-01  & 1.6782e+00   \\
$M=8 $ & 1.0461e-01  & 8.9300e-02  & 4.4795e-01   \\
$M=16$ & 2.5907e-02  & 2.2832e-02  & 1.1670e-01   \\
\hline   
order  & 2.023 & 1.942  & 1.923   \\
\hline
\hline
\end{tabular}
\end{center}
\end{table}
}
\end{example}

\section{Conclusions}
\label{sec-conclusion}
\setcounter{equation}{0}
We analyzed a linearized backward Euler Galerkin-mixed FEM for
the time-dependent Ginzburg--Landau equations
in both two and three dimensional spaces.
We have established unconditionally optimal error estimates 
for the three dimensional model, 
while analysis presented in this paper can be extended to many other cases, 
such as the two dimensional problem and Crank--Nicolson scheme. 
The method can solve for the induced magnetic field $\C \, \bA$ 
($\curl \, \bA$ in two dimensional space)
accurately without corner singularities and
the method is particularly suitable for nonconvex domains and 
domains with complex geometries.
Numerical experiments presented in this paper show the efficiency of 
the method and confirm our theoretical analysis. 
Large scale parallel computations and adaptive local refinement or 
moving mesh strategy will be conducted in our future work.

\section*{Acknowledgment}
The first author would like to thank Professor Douglas Arnold
for useful discussions on mixed methods for
vector Poisson equations.


\newpage

\section*{Appendix}
\setcounter{equation}{0}

With all the notations in section \ref{femmethod},
for any given $e_{\bA}^n \in \accentset{\circ}{\mathbf{U}}^{r}_{h}$,
if there exists a $e_{\bm{\sigma}}^n \in \accentset{\circ}{\mathbf{Q}}^{r+1}_{h}$ such that
\begin{eqnarray} 
  && (e_{\bm{\sigma}}^n, \bm{\chi}_h) - (\C \, \bm{\chi}_h, e_{\bA}^n) = 0, 
  \quad \forall  \, 
  \bm{\chi}_h \in \accentset{\circ}{\mathbf{Q}}^{r+1}_{h}\, ,
  \label{property-dis-curl}
\end{eqnarray}
then the following discrete Sobolev embedding inequality holds
\begin{eqnarray}
  \|e_{\bA}^{n}\|_{L^3} 
  & \le &
  C \left(\|e_{\bm{\sigma}}^{n}\|_{L^2} + \|\D \, e_\bA^{n}\|_{L^2} \right) \,.
  \label{embedding-results}
\end{eqnarray}

\noindent {\bf{ \em Proof:}}

From the Hodge decomposition \cite{AFW,AFW1},
the following exact sequence holds:
\begin{displaymath}
\begin{array}{ccccc}
  \accentset{\circ}{\bH}(\C)&  \xrightarrow{\,\, \C \,\,} & \accentset{\circ}{\bH}(\D)
  & \xrightarrow{\,\, \D \,\,}  & L_0^2       \\
  \downarrow&             &\downarrow &             &\downarrow  \\
  \accentset{\circ}{\mathbf{Q}}_h^{r+1} & \xrightarrow{\,\, \C \,\,}
  & \accentset{\circ}{\mathbf{U}}_h^r & \xrightarrow{\,\, \D \,\,}
  & \accentset{\circ}{\mathbf{S}}_h^r
\end{array}
\end{displaymath}
where $L_0^2 = \{v \in L^2 \, | \, (v,1)=0\}$
and $\accentset{\circ}{\mathbf{S}}_h^r={\mathbf{S}}_h^r \cap L_0^2$ with
\[
  {\mathbf{S}}_h^r:= \{v \in L^2(\Omega) \, | \,
  \textrm{$\forall \, \Omega_{K} \in \mathcal{T}_{h}$,
  $v|_{\Omega_{K}}$ is a polynomial of degree $r$}\}.
\]
Then, we have
\begin{eqnarray}
  \accentset{\circ}{\mathbf{U}}_h^r = \C \, \accentset{\circ}{\mathbf{Q}}_h^{r+1}
  \oplus \mathrm{ \bf grad}_h \accentset{\circ}{\mathbf{S}}_h^r \,,
  \label{hodge-discrete}
\end{eqnarray}
where $\oplus$ denotes the direct sum 
and the linear operator
$\mathrm{ \bf grad}_h:\accentset{\circ}{\mathbf{S}}_h^r \rightarrow
\accentset{\circ}{\mathbf{U}}_h^{r}$ is defined as:
for any given $s_h \in \accentset{\circ}{\mathbf{S}}_h^r$,
find $\mathrm{ \bf grad}_h s_h \in \accentset{\circ}{\mathbf{U}}_h^{r}$ such that
\begin{align}
  (\mathrm{ \bf grad}_h s_h, \bm{v}_h) = -(s_h, \D \, \bm{v}_h),
  \quad \forall \bm{v}_h \in \accentset{\circ}{\mathbf{U}}_h^{r} \,.
  \label{define-grad-h}
\end{align}
Therefore, for any given $e_{\bA}^{n} \in \accentset{\circ}{\mathbf{U}}_h^{r}$,
there exist $ \bm{\theta}_h \in \accentset{\circ}{\mathbf{Q}}_h^{r+1} $
and $s_h \in \accentset{\circ}{\mathbf{S}}_h^r$ such that
\begin{eqnarray}
  e_{\bA}^{n} = \C \, \bm{\theta}_h \oplus \mathrm{ \bf grad}_h s_h
  \label{decomposition-rt}
\end{eqnarray}
To prove \refe{embedding-results}, 
we only need to  show
\[
\|\mathrm{ \bf grad}_h s_h \|_{L^3} + \|\C \,\bm{\theta}_h \|_{L^3}
 \le 
C \left(\|e_{\bm{\sigma}}^{n}\|_{L^2} + \|\D \, e_\bA^{n}\|_{L^2} \right) \,.
\]

We first estimate $\mathrm{ \bf grad}_h s_h$.
Indeed, $\mathrm{ \bf grad}_h s_h$ can be viewed as the
mixed FEM solution to the following
Poisson equation with pure Neumann boundary condition
\begin{eqnarray}
  \left\{
  \begin{array}{ll}
    - \Delta \, u = \D \,e_{\bA}^{n}   & \textrm{in $\Omega$}\\
    \nabla u \cdot \mathbf{n} =  0 & \textrm{on $\partial \Omega$}
  \end{array}
  \right.
  \label{poisson-appendix}
\end{eqnarray}
It should be noted that the regularity for
\refe{poisson-appendix} depends upon the the domain $\Omega$,
see \cite{Dauge1,Dauge2}.
For a general polyhedron in three dimensional space,
we have the following shift theorem
\begin{eqnarray}
\| u \|_{H^{3/2+\alpha}} \le C \|\D \,e_{\bA}^{n}\|_{H^{-1/2+\alpha}} \,,
\label{poisson-shift}
\end{eqnarray}
where $\alpha>0$ depends only upon the geometry of the polyhedron.
The classical mixed FEM for solving  \refe{poisson-appendix} is to
find $(\bm{q}_h,u_h) \in
(\accentset{\circ}{\mathbf{U}}_h^r,\accentset{\circ}{\mathbf{S}}_h^r)$
such that
\begin{eqnarray}
  \left\{
  \begin{array}{ll}
    (\bm{q}_h, \bm{v}_h) + (u_h \, , \D \,\bm{v}_h)= 0
    & \forall  \, \bm{v}_h \in \accentset{\circ}{\mathbf{U}}_h^r \,, \\
    -(\D \,\bm{q}_h, \mu_h)  =  (\D \, e_{\bA}^{n}, \mu_h)
    & \forall \, \mu_h \in \accentset{\circ}{\mathbf{S}}_h^r \,.
  \end{array}
  \right.
  \nn
\end{eqnarray}
From  \refe{define-grad-h} and \refe{decomposition-rt},
one can verify that  $\mathrm{ \bf grad}_h s_h = \bm{q}_h$.
By using an inverse inequality and  the classical error estimates
of mixed methods for elliptic equation \cite{Braess,BF},
we can deduce that
\begin{eqnarray}
  \|\mathrm{ \bf grad}_h s_h \|_{L^3}
  & \le & \|\mathrm{ \bf grad}_h s_h - \pi_h\nabla u\|_{L^3}
  + \|\pi_h \nabla u\|_{L^3}
  \nn \\
  & \le & Ch^{-\frac{1}{2}}\|\mathrm{ \bf grad}_h s_h - \pi_h\nabla u\|_{L^2}
  + \|u\|_{H^{{3}/{2} + \alpha}}
  \nn \\
  & \le & Ch^{-\frac{1}{2}}
  (\|\mathrm{ \bf grad}_h s_h - \nabla u\|_{L^2} + \| \pi_h\nabla u - \nabla u\|_{L^2})
  + C\|\D \, e_{\bA}^{n}\|_{L^2}
  \nn \\
  & \le & Ch^{-\frac{1}{2}} h^{\frac{1}{2}+\alpha} {\left\|\nabla u \right\|}_{H^{{1}/{2}+\alpha}}
  + C\|\D \, e_{\bA}^{n}\|_{L^2}
  \nn \\
  & \le & C\|\D \, e_{\bA}^{n}\|_{L^2} \,,
  \label{estimate-1}
\end{eqnarray}
where $\pi_h$ is the projection operator defined in \refe{interpolate}.

Now we turn to estimate $\C \, \bm{\theta}_h$.
From the discrete Hodge decomposition \refe{decomposition-rt},
we have 
\begin{eqnarray}
  (\C \, \bm{\theta}_h , \C \,\bm{\chi}_h)
  =  (e_{\bA}^{n}, \C \, \bm{\chi}_h)
  = (e_{\bm{\sigma}}^{n}, \bm{\chi}_h)\,
  \quad \textrm{for any $\bm{\chi}_h \in \accentset{\circ}{\mathbf{Q}}_h^{r+1}$.}
  \nn
\end{eqnarray}
A key observation is that
$\bm{\theta}_h$ can be viewed as the mixed finite element solution to
the vector Poisson equation
\begin{eqnarray}
  \left\{
  \begin{array}{ll}
    \C \, \C \, \bm{u} - \nabla \, \D \, \bm{u} =  e_{\bm{\sigma}}^{n} & \textrm{in $\Omega$}\\
    \bm{u} \times \mathbf{n} =  \mathbf{0}, \D \bm{u} = 0, & \textrm{on $\partial \Omega$}
  \end{array}
  \right.
  \label{vector-Poisson}
\end{eqnarray}
which is to find
$(\phi_h,\bm{u}_h) \in \accentset{\circ}{V}_h^{r}\times \accentset{\circ}{\mathbf{Q}}_h^{r+1}$,
where $\accentset{\circ}{V}_h^{r}$ denote the $r$-th Lagrange element space
with zero trace, such that
\begin{eqnarray}
  \left\{
  \begin{array}{ll}
    -( \phi_h, \zeta_h) + ( \bm{u}_h \, \nabla \zeta_h) = 0 ,&
    \textrm{for any $\zeta_h \in \accentset{\circ}{V}_h^{r}$.} \\
    (\nabla \phi_h, \bm{\chi}_h)
    + (\C \, \bm{u}_h , \C \,\bm{\chi}_h) =(e_{\bm{\sigma}}^{n}, \bm{\chi}_h)\,,
    & \textrm{for any $\bm{\chi}_h \in \accentset{\circ}{\mathbf{Q}}_h^{r+1}$.}
  \end{array}
  \right.
  \label{vector-Poisson-fem}
\end{eqnarray}
Since $\bm{\theta}_h = \bm{u}_h$, from the standard error estimate \cite{AFW, AFW1},
we have
\begin{align}
& \|\D \bm{u} - \phi_h\|_{L^2} 
+ \|\bm{u} - \bm{\theta}_h\|_{\mathbf{H}(\C)} 
\nn \\
& \le C \Big(
\inf_{\zeta_h \in \accentset{\circ}{V}_h^{r}} \| \D \bm{u} - \zeta_h \|_{L^2}
+\inf_{\bm{\chi}_h \in \accentset{\circ}{\mathbf{U}}_h^r}
\|\mathbf{u} - \bm{\chi}_h\|_{\mathbf{H}(\C)} 
\Big)
\end{align}
Thus, with the projection operator $\pi_h$ in \refe{interpolate},
we can derive the following estimates
\begin{eqnarray}
  \|\C \,\bm{\theta}_h \|_{L^3}
  & \le & \|\C \,\bm{\theta}_h - \C \, \pi_h\bm{u} \|_{L^3} + \|\C \, \pi_h\bm{u}\|_{L^3}
  \nn \\
  & \le &  C h^{-\frac{1}{2}} \|\C \,\bm{\theta}_h - \C \, \bm{u}\|_{L^2}
  +  C h^{-\frac{1}{2}}\| \C \, \pi_h\bm{u} - \C \, \bm{u}\|_{L^2}
  + \|e_{\bm{\sigma}}^{n}\|_{L^2}
  \nn \\
  & \le & C h^{-\frac{1}{2}}\| \C \, \pi_h\bm{u} - \C \, \bm{u}\|_{L^2}
  + \|e_{\bm{\sigma}}^{n}\|_{L^2}
  \nn \\
  & \le & C h^{-\frac{1}{2}} h^{\frac{1}{2}+\alpha}
  (\|\C \, \bm{u}\|_{H^{{1}/{2}+\alpha}} + \|\D \, \bm{u}\|_{H^{{1}/{2}+\alpha}})
  + \|e_{\bm{\sigma}}^{n}\|_{L^2}
  \nn \\
  & \le & C \| e_{\bm{\sigma}}^{n} \|_{L^2} \,,
  \label{estimate-2}
\end{eqnarray}
where we have used an inverse inequality and noted that in \refe{vector-Poisson-fem}
\begin{eqnarray}
  \|\C \, \bm{u}\|_{H^{{1}/{2}+\alpha}} + \|\D \, \bm{u}\|_{H^{{1}/{2}+\alpha}}
  \le C \| e_{\bm{\sigma}}^{n} \|_{L^2}
\label{embedding-continuous}
\end{eqnarray}
for a certain $\alpha>0$, see \cite{Alonso_Valli,ABDG}.

The embedding inequality \refe{embedding-results} is proved by
combining  \refe{estimate-1} and \refe{estimate-2}. \qed

\begin{remark}

  The discrete embedding inequality \refe{embedding-results} is enough in
  the error analysis of the Galerkin-mixed FEM for the TDGL equations.
  However, if the polyhedron $\Omega$ is convex
  we can derive a stronger discrete embedding inequality 
  \begin{eqnarray}
  \|e_{\bA}^{n}\|_{L^6} 
  & \le &
  C \left(\|e_{\bm{\sigma}}^{n}\|_{L^2} + \|\D \, e_\bA^{n}\|_{L^2} \right) \,.
  \end{eqnarray}
  On a convex domain,  we have $\alpha > \frac{1}{2}$ in \refe{poisson-shift}
  for the shift theorem. And then, we have
  \begin{eqnarray}
    \|\mathrm{ \bf grad}_h s_h \|_{L^6}
    & \le & \|\mathrm{ \bf grad}_h s_h - \pi_h\nabla u\|_{L^6}
    + \|\pi_h \nabla u\|_{L^6}
    \nn \\
    & \le & Ch^{-1}\|\mathrm{ \bf grad}_h s_h - \pi_h\nabla u\|_{L^2}
    + \|u\|_{H^{2}}
    \nn \\
    & \le & C\|\D \, e_{\bA}^{n}\|_{L^2} \,.
  \end{eqnarray}
  Moreover, for a convex domain $\Omega$,
  the embedding inequality \refe{embedding-continuous} can be re-written by
  \begin{eqnarray}
    \|\C \, \bm{u}\|_{H^1} + \|\D \, \bm{u}\|_{H^1}
    \le C \| e_{\bm{\sigma}}^{n} \|_{L^2},
    \label{embedding-continuous-full}
  \end{eqnarray}
  which implies   
  \begin{eqnarray}
    \|\C \,\bm{\theta}_h \|_{L^6} \le  C \| e_{\bm{\sigma}}^{n} \|_{L^2}.
  \end{eqnarray}
\end{remark}

\begin{remark}

  The discrete Hodge decomposition \refe{hodge-discrete}
  is only valid for simply-connected domain.
  If the domain $\Omega$ is multi-connected, we have
  \begin{eqnarray}
    \accentset{\circ}{\mathbf{U}}_h^r = \C \, \accentset{\circ}{\mathbf{Q}}_h^{r+1}
    \oplus \mathrm{ \bf grad}_h \accentset{\circ}{\mathbf{S}}_h^r
    \oplus \mathrm{ \bf Har}_h^r \,,
  \end{eqnarray}
  where
  $\mathrm{\bf Har}_h^r := \{ \bm{v}_h \in \accentset{\circ}{\mathbf{U}}_h^r \,
  | \, \D \bm{v}_h = 0, (\bm{v}_h \, , \C \bm{\chi}_h) = 0, \forall \bm{\chi}_h
  \in \accentset{\circ}{\mathbf{Q}}_h^{r+1}\}$
  stands for the harmonic part of
  $\accentset{\circ}{\mathbf{U}}_h^r$. In this case,
  \begin{eqnarray}
    e_{\bA}^{n} = \C \, \bm{\theta}_h
    \oplus \mathrm{ \bf grad}_h s_h
    \oplus \bm{v}_h \,,
    \nn
  \end{eqnarray}
  where $\bm{v}_h \in \mathrm{\bf Har}_h^r$.
  we can prove the following inequality via a similar analysis
  \[
  \|\bm{v}_h\|_{L^3}
  \le C \left(\|e_{\bm{\sigma}}^{n}\|_{L^2} + \|\D \, e_\bA^{n}\|_{L^2} \right) \,.
  \]
\end{remark}

\end{document}